\documentclass[11pt]{article}
\usepackage{IEEEtrantools}
\usepackage{mathtools, nccmath}
\usepackage{latexsym}
\usepackage{amsfonts}
\usepackage{amssymb}
\usepackage{eqnarray, amsmath,amsfonts,amsthm}
\usepackage{mathrsfs}
\usepackage{enumerate}
\usepackage{dsfont}
\usepackage{hyperref}
\usepackage{graphicx}
\usepackage{epstopdf}
\usepackage{tikz}

\usetikzlibrary{shapes,arrows.meta,positioning,calc}
\usepackage{cleveref}
\usepackage{tikz-cd}
\usepackage{psfrag}

\usepackage{multirow}
\usepackage{tikz}
\usetikzlibrary{arrows,decorations.markings}
\usepackage{bookmark}
\usepackage{hyperref}
\hypersetup{
     colorlinks   = true,
     citecolor    = blue
}

\newcommand{\newsection}[1]{\setcounter{equation}{0}
\setcounter{dfn}{0}
\section{#1}}

\newtheorem{dfn}{Definition}[section]
\newtheorem{thm}[dfn]{Theorem}
\newtheorem{lmma}[dfn]{Lemma}
\newtheorem{ppsn}[dfn]{Proposition}

\newtheorem{rmrk}[dfn]{Remark}
\newtheorem{notation}[dfn]{Notation}
\newtheorem*{theoremA}{Theorem A}
\newtheorem*{theoremB}{Theorem B}
\newtheorem*{theoremC}{Theorem C}
\newtheorem*{theoremD}{Theorem D}

\newcommand{\ot}{\otimes}

\newcommand{\dt}{\Delta_n}
\newcommand{\bbc}{\mathbb{C}}

\newcommand{\bbn}{\mathbb{N}}

\def \qed { \mbox{}\hfill
$\Box$\vspace{1ex}}

\title{Sections and Chapters}
\begin{document}

\tikzset{->-/.style={decoration={
  markings,
  mark=at position #1 with {\arrow{>}}},postaction={decorate}}}
  \tikzset{-<-/.style={decoration={
  markings,
  mark=at position #1 with {\arrow{<}}},postaction={decorate}}}

\author{\sc{Keshab Chandra Bakshi\footnote{\, Partially supported by DST INSPIRE Faculty grant DST/INSPIRE/04/2019/002754}\,,\,\,Satyajit Guin,\,\,Guruprasad}}
\title{A noncommutative construction of families of biunitary matrices and application to subfactors}
\maketitle


\begin{abstract}
We introduce a construction that, given a pair $(u,v)$ of complex Hadamard matrices of the same order, generates infinitely many {\it biunitary} matrices of varying (and distinct) orders. As a key application, this framework yields nested sequences of vertex model subfactors that are not a tower of downward basic construction. Notably, the construction is noncommutative: interchanging the matrices (i.e., considering $(v,u)$ instead of $(u,v)$) can lead to non-isomorphic subfactors. Focusing on the Hadamard equivalence class of the Fourier matrix $F_n$, we provide a full characterization of the resulting vertex model subfactors, along with explicit computations of their relative commutants. Along the way, we conduct a detailed study of certain naturally arising inner and outer automorphisms that play a key role in the structure of these subfactors.
\end{abstract}
\bigskip

{\bf AMS Subject Classification No.:} {\large 46}L{\large 37}, {\large 46}L{\large 55}
\smallskip


\hypersetup{linkcolor=blue}
\tableofcontents
\bigskip


\newsection{Introduction}\label{Sec1}

Since the early development of subfactor theory \cite{Jo}, matrices (more precisely, multi-matrices) have played a central role in the construction of subfactors—particularly hyperfinite ones. The general strategy involves starting with a matrix as initial data and using it to produce a non-degenerate commuting square of finite-dimensional  $C^*$-algebras. Applying Jones' basic construction to this square then yields a hyperfinite subfactor. This naturally leads to a fundamental question: how much of the subfactor’s structure can be recovered or understood from the input matrix? In this article, we focus on two notable classes where this connection is especially rich and relevant. The first is the case of {\em Hadamard} (or {\em spin model}) subfactors. Any complex Hadamard matrix $u\in M_n(\bbc)$ induces the following non-degenrate commuting square
\[
\begin{matrix}
\mathrm{Ad}_u(\Delta_n) &\subset & M_n(\bbc)\\
\cup & & \cup\\
\bbc & \subset & \Delta_n
\end{matrix}
\]
where $\Delta_n$ denotes the subalgebra of diagonal matrices in $M_n(\bbc)$, and then, iterating Jones' basic construction (vertically), one gets a special type of hyperfinite subfactor $R_u\subset R$, known as the Hadamard/spin model subfactor. This class of subfactors are always irreducible with integer Jones index (in fact, $[R:R_u]=n$). Their importance has been emphasized by Jones on several occasions (see \cite{JS, Jo2}, for instance). For recent developments on spin model subfactors, the readers are referred to \cite{BG1, BG2, BGG, Bur, M, N}.

The second important class of subfactors that we deal with in this article is the {\it vertex model} subfactors. Recall that a unitary matrix $w=(w_{\alpha a}^{\beta b})\in M_n(\bbc)\otimes M_k(\bbc)$ is called a {\it biunitary} matrix if its block transpose $\widetilde{w}= (\widetilde{w}_{\alpha a}^{\beta b})$, with $\widetilde{w}_{\alpha a}^{\beta b}:=w_{\beta a}^{\alpha b}$, is also a unitary matrix in $M_n(\bbc)\otimes M_k(\bbc)$ (see \cite{BINA,JS}). A biunitary matrix $w\in M_n(\bbc)\otimes M_k(\bbc)$ induces the following non-degenerate commuting square
\[
\begin{matrix}
\mathrm{Ad}_w(M_n(\bbc)\otimes\bbc) &\subset & M_n(\bbc)\otimes M_k(\bbc)\\
\cup & & \cup\\
\bbc & \subset & \bbc\otimes M_k(\bbc)
\end{matrix}
\]
and then, iterating Jones' basic construction (vertically), one gets another special type of hyperfinite subfactor $R_w\subset R$, known as the vertex model subfactor. This class of subfactors are not necessarily irreducible, unlike the spin model subfactors, and these are of integer Jones index (in fact, $[R:R_w]=k^2$). For earlier works on vertex model subfactors, see \cite{Ba, BINA, KSV, KS, R}.

Both of these types of subfactors are of great significance; however, relatively little is known about their detailed structure. Inspired by Jones's two subfactor theory \cite{Jo1}, the investigation into pairs of spin model subfactors is initiated in \cite{BG1, BG2, BGG}. This article serves as a continuation of that line of inquiry, presenting new insights. Specifically, we establish a curious relationship between these two classes of subfactors in that we construct infinitely many non-isomorphic subfactors of the second kind from a fixed pair (but with varying Jones index) of the first kind. Indeed, we start with two complex Hadamard matrices of the same order and obtain infinitely many nondegenerate commuting squares of finite-dimensional $C^*$-algebras induced by a family of biunitary matrices. Performing Jones' basic construction, we then get a nested sequence of infinitely many non-isomorphic vertex model subfactors. Observe that two complex Hadamard matrices give two non-degenerate commuting squares, and so, naturally produce two spin model subfactors. In \cite{BG2}, using two spin model subfactors each of index $4$, the first two authors constructed infinitely many subfactors of $R$, many of which are of infinite depth (see \cite{BH} in this regard). The ability to generate infinitely many subfactors from just a pair of matrices is, in itself, a surprising phenomenon—a theme that carries over into the present work. Our first result is the following (see Theorem $3.2$, \Cref{Sec3}).
\begin{theoremA}
Let $n\in\bbn$ and $(u,v)$ be a pair of complex Hadamard matrices of order $n$. Then, there is a sequence $\{BU(u,v;k)\}_{k=0}^\infty$ of biunitary matrices such that the order of $BU(u,v;k)$ is $n^{k+2}$.
\end{theoremA}
Traversing from finite to infinite dimension, as an application of Theorem A, we prove the following result (see Theorem $4.1$, \Cref{Sec4}).
\begin{theoremB}
Let $n\in\bbn$ and $u$ be a complex Hadamard matrix of order $n$. Let $R_u\subset R$ be the corresponding Hadamard/spin model subfactor of the hyperfinite type $II_1$ factor $R$. Choose another complex Hadamard matrix $v$ of the same order $n$. Then, for each $k\in\mathbb{N}\cup\{0\}$, there exists a vertex model subfactor $R^{u,v}_k \subset R$ such that
\begin{enumerate}[$(i)$]
\item $R^{u,v}_k\subset R_u$ for each $k$, and the Jones index is given by $[R:R^{u,v}_k]=n^{2(k+1)}$;
\item The subfactors $\{R^{u,v}_k\subset R:k\in\mathbb{N}\cup\{0\}\}$ satisfy the inclusion $R^{u,v}_{k+1} \subset R^{u,v}_k$ for each $k$.
\end{enumerate}
Thus, we have a nested sequence of vertex model subfactors
\[
\cdots \subset R^{u,v}_{k+1}\subset R^{u,v}_k\subset\cdots\subset R^{u,v}_0\subset R
\]
each having the Hadamard/spin model subfactor $R_u\subset R$ as an intermediate subfactor.
\end{theoremB}
We show that the nested sequence of subfactors $\{R^{u,v}_k\subset R\}_{k=1}^\infty$ is not the downward basic construction of the first vertex model subfactor $R^{u,v}_0\subset R$, and the first vertex model subfactor $R^{u,v}_0\subset R$ itself is not the downward basic construction of the Hadamard subfactor $R_u\subset R$ and so on. Therefore, our construction is new. Moreover, our construction is noncommutative in nature; if we switch the position of $u\mbox{ and }v$, one can obtain non-isomorphic subfactors $(R^{u,v}_k\subset R) \ncong (R^{v,u}_k\subset R)$. This also says that the biunitary matrices $BU(u,v;k)$ in Theorem A need not be conjugate to $BU(v,u;k)$ (otherwise, the resulting subfactors would have been always isomorphic). We would like to highlight one crucial point here. Upon obtaining Theorem A, one might assume that infinitely many vertex model subfactors can be immediately constructed via the iterated basic construction. However, the key contribution of Theorem B is that it allows us to construct a {\em nested sequence} of vertex model subfactors that is not a tower of downward basic construction, all of which are contained within the hyperfinite factor $R$.

Both Theorems A and B can be most effectively illustrated diagrammatically as follows\,:
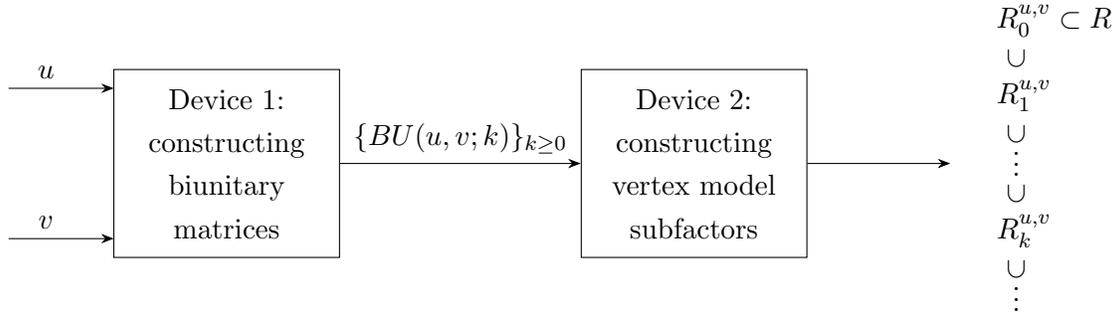
\begin{center}
\begin{tikzpicture}[>=Stealth, font=\normalsize]
\node[draw, minimum width=3cm, minimum height=2.5cm, align=center] (device1) at (3.7,0) {Device 1:\\ constructing\\ biunitary\\ matrices};
\node[draw, minimum width=3cm, minimum height=2.5cm, align=center, right=3.2cm of device1] (device2) {Device 2:\\ constructing\\ vertex model\\ subfactors};
\draw[->] (.8,1) -- (2.2,1);
\node[above] at (1.3,1) {$u$};
\draw[->] (.8,-1) -- (2.2,-1);
\node[above] at (1.3,-1) {$v$};
\draw[->] (device1.east) -- (device2.west);
\node[above] at ($(device1.east)!0.5!(device2.west)$) {$\{BU(u,v;k)\}_{k\geq 0}$};
\coordinate (out0) at ([xshift=1.9cm]device2.east);
\draw[->] (device2.east) -- (out0);
\node at (14.7,1.9)  {{$R^{u,v}_0 \subset R$}};
\node at (14.2, 1.4) {$\rotatebox{90}{$\subset$}$};
\node at (14.3, 0.9)  {{$R^{u,v}_1$}};
\node at (14.2,0.4) {$\rotatebox{90}{$\subset$}$};
\node at (14.2, 0.1)  {{$\vdots$}};
\node at (14.2,-0.4) {$\rotatebox{90}{$\subset$}$};
\node at (14.3, -0.9)  {{$R^{u,v}_k$}};
\node at (14.2,-1.4) {$\rotatebox{90}{$\subset$}$};
\node at (14.2,-1.7)  {{$\vdots$}};
\end{tikzpicture}
\end{center}

\noindent Device $1$ is precisely the Theorem A that takes input as $(u,v)$ and produces biunitary matrices $\{BU(u,v;k)\}_{k\geq 0}$. Then, we enter into Device $2$ that manufactures subfactors from matrices, and as final output, we get the nested sequence of subfactors $\cdots\subset R^{u,v}_k\subset\cdots\subset R$ which is precisely Theorem B.

If we consider two such pairs $(u,v),\,(\widetilde{u},\widetilde{v})$ as our input data, then our device produces infinitely many pairs of vertex model subfactors $R_k^{u,v},R_k^{\widetilde{u},\widetilde{v}}\subset R$ for each $k\in\bbn\cup\{0\}$. A natural question that immediately arises in this context is the following\,:
\smallskip

\noindent\textbf{Question\,:} Are the pairs of vertex model subfactors $R_k^{u,v},R_k^{\widetilde{u},\widetilde{v}}\subset R$ isomorphic for all $k\in\bbn\cup\{0\}$?
\smallskip

\noindent This question is investigated in \Cref{Sec5} and we have the following result (See Theorem $5.2$, \Cref{Sec5}).
\begin{theoremC}
Let $(u,v),(\widetilde{u},\widetilde{v})$ be pairs of complex Hadamard matrices all of the same order. Consider the pairs of vertex model subfactors $R_k^{u,v},\,R_k^{\widetilde{u},\widetilde{v}}\subset R,\,k\in\bbn\cup\{0\},$ constructed in Theorem B with the pairs of Hadamard subfactors $R_u,R_{\widetilde{u}}\subset R$ as their intermediate subfactors respectively. Then,
\begin{enumerate}[$(i)$]
\item If $v=\widetilde{v}$ and $u\simeq\widetilde{u}$, where $\simeq$ denotes the Hadamard equivalence, then there exists $\theta\in\mathrm{Aut}(R)$ such that $\big(R_k^{u,v}\subset R\big)\cong\big(R_k^{\widetilde{u},\widetilde{v}}\subset R\big)$ via $\theta$ for each $k\in\bbn\cup\{0\}$.
\item If $v\neq\widetilde{v}$ or $u\not\simeq\widetilde{u}$, the pairs of vertex model subfactors need not be isomorphic.
\end{enumerate}
\end{theoremC}
A striking feature of part $(i)$ in Theorem C is that one single automorphism $\theta$ works for all the subfactors; in other words, $\theta$ implements a chain isomorphism. Using Connes' machinery from \cite{Connes}, we are also able to completely characterize when the automorphism $\theta$ is inner and outer, along with its outer period.

In the final section of the paper, taking help of Theorem C, we provide a complete characterization of the vertex model subfactors constructed in Theorem B arising from the equivalence classes of the DFT (Discrete Fourier Transform) matrices. It turns out that all such subfactors are diagonal. Furthermore, the subfactors derived from the equivalence classes of the DFT matrices exhibit a certain ``rigidity" when we restrict our attention to diagonal subfactors. More precisely, once we move beyond DFT matrices, the resulting vertex model subfactors need not be diagonal. We have the following result (see Theorem $6.1$, \Cref{Sec6}).
\begin{theoremD}
Let $F_n$ be the DFT matrix of order $n$ and $[F_n]$ denote the Hadamard equivalence class of $F_n$. We have the following\,:
\begin{enumerate}[$(i)$]
\item If $u\in[F_n]$, then there exists $\theta\in\mathrm{Out}(R)$ with outer period $n$ such that the vertex model subfactor $R_k^{u,F_n}\subset R$ obtained in Theorem B is isomorphic to the following diagonal subfactor
\[
\mathscr{R}_k:=\bigg{\{} \displaystyle{\sum_{j=1}^{n^{k+1}} E_{jj} \ot \theta^{\,r_j}(x) : \, x \in R,\,r_j\in\{0,\ldots,n-1\}},\,\theta^n=\mathrm{id}\bigg{\}} \subset M^{(k+1)}_n\ot R
\]
for any $k\geq 0$. Furthermore, the relative commutant is given by the following\,:
\[
\big(R^{u, F_n}_k\big)^\prime\cap R=\bbc^{n^{k+1}}.
\]
\item If $v\neq F_n$ or $u\notin[F_n]$, the vertex model subfactors $R_k^{u,v}\subset R$ need not be diagonal.
\end{enumerate}
\end{theoremD}
Again, a striking feature of part $(i)$ in Theorem D is that one single outer automorphism works for the entire chain of the subfactors. This result indicates that our construction is fully understood in the case of the DFT matrices. That the subfactors are not always diagonal (part $(ii)$, Theorem C) if $u\notin[F_n]$ is, in fact, easy to see. Choose $u$ such that the Hadamard subfactor $R_u\subset R$ is of infinite depth (see \cite{Bur} for instances of such). Since $R_u\subset R$ is an intermediate subfactor for each vertex model subfactor $R_k^{u,v}\subset R$, all these are must be of infinite depth by \cite{Bi}, irrespective of the choice of $v$. Therefore, these cannot be diagonal subfactors. However, the case of $v\neq F_n$ in part $(ii)$ requires some work.
\medskip

Throughout the article, the underlying field is always $\bbc$, and hence, we omit it for notational brevity. Also, the complex number `$\textbf{i}$' shall be explicitly denoted by $\sqrt{-1}$. This is to avoid any possible confusion with the index `$i$' frequently used in the article.


\newsection{Preliminaries}\label{Sec2}

In this section, we review some key concepts used throughout the paper. While we aim to make the article as self-contained as possible, we avoid delving too deeply into these concepts. However, relevant citations are provided at appropriate places for the readers who wish to explore the material further.

\subsection{Complex Hadamard matrices and biunitary matrices}\label{Sec2.1}

Two special types of matrices are relevant in this article, namely, the {\it complex Hadamard} matrices and the {\it biunitary} matrices. Let us begin with the first one.

Recall that a matrix with entries $\pm 1$ and mutually orthogonal rows and columns is called a Hadamard matrix. If $H$ is a Hadamard matrix and $H^{\intercal}$ denote its transpose, then $HH^{\intercal}=n I_n$. A complex Hadamard matrix is a generalization of a Hadamard matrix over the complex field. A complex Hadamard matrix is a matrix in $M_n(\bbc)$, henceforth simply denoted by $M_n$ throughout the article, with all the entries of the same modulus and $HH^*=nI_n$. Observe that $\frac{1}{\sqrt{n}}H$ is a unitary matrix. In this article, we work with the following alternative definition of the complex Hadamard matrices.
\begin{dfn}[\cite{JS}]
A complex Hadamard matrix of order $n$ is a unitary matrix in $M_n$ such that each of its entries has modulus $\frac{1}{\sqrt{n}}$.
\end{dfn}
We denote by $\mathbb{H}(M_n)$ the set of all complex Hadamard matrices in $M_n$. For every $n\geq 1$, the Fourier matrix $(F_n)_{i,j}:=\frac{1}{\sqrt{n}}\exp\big(2\pi\sqrt{-1}(i-1)(j-1)/n\big),\,1\leq i,j\leq n,$ is a complex Hadamard matrix, widely referred as a DFT (discrete Fourier transform) matrix. Two complex Hadamard matrices are called Hadamard equivalent, denoted by ${\displaystyle H_{1}\simeq H_{2}}$, if there exist diagonal unitary matrices ${\displaystyle D_{1}, D_{2}}$ and permutation matrices ${\displaystyle P_{1}, P_{2}}$ such that
\begin{eqnarray}\label{Hadamard equivalence}
H_{1}=D_{1}P_{1}H_{2}P_{2}D_{2}\,.
\end{eqnarray}
This is an equivalence relation on the set $\mathbb{H}(M_n)$. For $n=2,3,5$, every complex Hadamard matrix is Hadamard equivalent to the Fourier matrix $F_n$; whereas, for $n=4$, it is known that there is a one-parameter continuous family of Hadamard inequivalent complex Hadamard matrices. Note that the complete classification of complex Hadamard matrices of higher orders ($n\geq 6$) is still unknown. 
\smallskip

The second class of matrices relevant to us is the so-called biunitary matrices.

\begin{dfn}[\cite{JS, BINA}]
A unitary matrix $u = (u_{\alpha a}^{\beta b}) \in M_n \otimes M_k$, is called a biunitary matrix if its block transpose $w = (w_{\alpha a}^{\beta b})$, with $w_{\alpha a}^{\beta b} := u_{\beta a}^{\alpha b}$, is also a unitary matrix in $M_n \otimes M_k$.
\end{dfn}
An effective tool that sometimes helps to avoid tedious computation to check whether a given unitary matrix in $M_n\otimes M_k$ is a biunitary matrix is in terms of the so-called `commuting square'. Recall that there is a unique normalized tracial state $tr$ on the type $I_{nk}$ factor $M_n\otimes M_k=M_{nk}$. Consider the unital subalgebras $\bbc\otimes M_k$ and $\mbox{Ad}_u(M_n\otimes\bbc)$ of $M_n\otimes M_k$, where $u\in\mathcal{U}(M_{nk})$ is a unitary matrix of order $nk$. Let $E_1$ and $E_2$ be the unique trace preserving conditional expectations onto $\bbc\otimes M_k$ and $\mbox{Ad}_u(M_n\otimes\bbc)$ respectively. Recall the following result from \cite{JS}.

\begin{thm}[\cite{JS}]\label{biunitary thm}
Let $u=(u_{\alpha a}^{\beta b})\in M_n\otimes M_k$ be a unitary matrix. Then, $u$ is a biunitary matrix if and only if the following quadruple
\begin{eqnarray}\label{biunitary matrix}
&&\begin{matrix}
\mathrm{Ad}_u(M_n\otimes\bbc) &\subset^{\,E_2} & M_n\otimes M_k\\
\cup^{\,\mbox{tr}} & &\cup^{\,E_1}\\
\bbc &\subset^{\,\mbox{tr}} & \bbc\otimes M_k
\end{matrix}
\end{eqnarray}
is a commuting square. That is, $\mbox{Ad}_u(M_n\otimes\bbc)\cap(\bbc\otimes M_k)=\bbc$ and $E_1E_2=E_2E_1=tr$, where `$tr$' is the unique normalized trace on $M_{nk}$.
\end{thm}

\subsection{Basic construction of finite-dimensional algebras}\label{Sec2.2}

\textbf{Finite-dimensional inclusion\,:} Let $B\subset A$ be a unital inclusion of finite-dimensional von Neumann algebras. Let
\[
B =\bigoplus_{i=1}^\ell M_{n_i}\hspace{2mm} \text{and} \hspace{2mm} A =\bigoplus_{j=1}^k M_{m_j},
\]
where $M_{n_i}$ (resp., $M_{m_j}$) denotes the type $I_{n_i}$ (resp., type $I_{m_j}$) factor.

\noindent \textbf{Inclusion Matrix\,:} The inclusion matrix $\Lambda_B^A$ (or simply $\Lambda$) is defined by the $k \times \ell$ matrix $\Lambda = (\Lambda_{ij})$, where $\Lambda_{ij}$ {\em essentially} represents the number of times the $i$-th summand of $B$ appears in the $j$-th summand of $A$.

\noindent \textbf{Markov Trace\,:} Since $A\cong\bigoplus_{j=1}^k M_{m_j}$, any trace on $A$ is uniquely determined by a vector $v = (v_1, v_2, \ldots, v_k)^T$, called the trace vector.
A trace $\tau$ on $A$ is called a \textit{Markov trace} for the unital inclusion $B\subset A$ if it satisfies the Markov condition:
\begin{equation*}
\Lambda^T\Lambda v=||\Lambda||^2v.
\end{equation*}
For other equivalent conditions, see \cite{JS}.

\noindent \textbf{Connected inclusion\,:} The inclusion $B\subset A$ is called connected if the Bratteli diagram for it is connected.

\noindent \textbf{Basic construction\,:} Suppose that $B\subset A$ is a unital inclusion of finite-dimensional unital von Neumann algebras equipped with a fixed faithful tracial state $\tau$ on $A$. The GNS Hilbert space of $A$ for the trace $\tau$ is naturally identified with $A$. The orthogonal projection $e_B$ from $A$ onto its closed subspace $B$ induces a unique trace $\tau$-preserving conditional expectation $E^A_B$ from $A$ onto $B$ given by $E^A_B(x)=y$, where $e_B(x)=y$. The unital von Neumann subalgebra $\langle A, e_B \rangle$, generated by $A$ and $e_B$, of $\mathcal{B}(A)$ is called Jones' basic construction for the unital inclusion $B\subset A$ with the Jones projection $e_B$. 

We denote any quadruple
\begin{eqnarray}\label{intial cs}
\begin{matrix}
A_{0} &\subset & A_{1} \cr \cup &\ & \cup\cr B_{0} &\subset &  \,B_{1}\,\,
\end{matrix}  
\end{eqnarray}
of finite-dimensional unital von Neumann algebras by $(B_{0} \subset B_{1}, \, A_{0} \subset A_{1})$ for notational brevity. If the quadruple is a non-degenerate (that is, $A_0B_1=B_1A_0=A_1$) commuting square and all the inclusions are connected, the unital inclusion $A_{0} \subset A_{1}$ has a unique Markov trace (Corollary 3.2.5, \cite{JS}). We have the finite-dimensional von Neumann algebra $A_{2}= \langle A_{1}, e_{2} \rangle $ obtained as the basic construction for the inclusion $A_{0} \subset A_{1}$ with Jones' projection $e_2$. The inclusion $B_{0} \subset  B_{1}  \subset B_{2}:=\langle B_{1} , e_2 \rangle$ is an instance of basic construction. We have unique Markov-trace, denoted by $\tau_{1}$, for the unital inclusion $A_{1} \subset A_{2}$. Further, the following quadruple
\begin{eqnarray*}
\begin{matrix}
 A_{1} &\subset & A_{2} \cr \cup &\ & \cup & & \cr  B_{1} &\subset &  \,B_{2}.\,\,
\end{matrix}  
\end{eqnarray*}
is also a non-degenerate commuting square with respect to the Markov trace $\tau_1$-preserving conditional expectations $E^{A_{2}}_{A_{1}}, \, E^{A_{2}}_{B_{2}}$ and $E^{A_{2}}_{B_{1}}$. Therefore, iterating the basic construction of the non-degenerate commuting square in \Cref{intial cs}, we obtain the following tower of non-degenerate commuting squares
\begin{eqnarray*}
\begin{matrix}
A_{0} & \subset & A_{1} &\subset & A_{2} & \subset & A_{3} & \cdots \cr \cup &\ & \cup & & \cup && \cup \cr  B_{0} &\subset &  \,B_{1} & \subset & B_{2} & \subset & B_{3} & \cdots  \,\,
\end{matrix}  
\end{eqnarray*}
of finite-dimensional unital von Neumann algebras (see $1.1.6$, \cite{Po2}, in this regard). We have a unique trace $tr$ on $\cup_{k} A_{k}$, and let $\mathcal{H}$ be the GNS-completion of $\cup_{k} A_{k}$ with respect to the trace $tr$. Then, $\overline{\bigcup_{k}A_k}^{\text{sot}},\,\overline{\bigcup_{k}B_k}^{\text{sot}}\subset\mathcal{B}(\mathcal{H})$ are type $II_1$-factors, and we have the following result.
\begin{ppsn}[\cite{JS}]\label{prop} 
Let $(B_{0} \subset B_{1}, \, A_{0} \subset A_{1})$ be a quadruple of a non-degenerate commuting square of finite-dimensional unital von Neumann algebras such that all the inclusions
are unital and connected. For the subfactor $N:=\overline{\bigcup_{k}B_k}^{\text{sot}}\subset M:=\overline{\bigcup_{k}A_k}^{\text{sot}}$, we have the following\,:
\begin{enumerate}[$(i)$]
\item $[M :\, N]=||\Lambda||^{2}$, where $\Lambda$ is the inclusion matrix for the inclusion $B_{0} \subset A_{0}$.
\item (Ocneanu compactness) \,$N^\prime\cap M=B_{1}^\prime\cap A_{0}$.
\end{enumerate}
\end{ppsn}

\subsection{Hadamard/spin model and vertex model subfactors}\label{Sec2.3}

The special types of matrices introduced in \Cref{Sec2.1} give rise to two distinct classes of subfactors\,: the Hadamard/spin model subfactors and the vertex model subfactors. We briefly recall them from \cite{JS}.

Let $\Delta_n$ denote the unital subalgebra of $M_n$ consisting of diagonal matrices. This is a {\it Masa} (maximal abelian self-adjoint subalgebra) in $M_n$. Any $n\times n$ complex Hadamard matrix $u$ induces the following non-degenerate commuting square (see \cite{JS})
\[
\begin{matrix}
 \Delta_n &\subset & M_n \cr \cup &\ & \cup\cr  \mathbb{C} &\subset &  u\Delta_n u^*\,\,.
\end{matrix}
\]
The basic construction of $\Delta_n\subset M_n$ is $\Delta_n\otimes M_n$ with the Jones' projection $e_2:=\sum_{i=1}^n E_{ii}\ot E_{ii}$, and that of $\Delta_n\otimes M_n$ is $M_n\otimes M_n$ with the Jones' projection $e_3:=J_n\otimes I_n$, where $J_n=\frac{1}{n}\sum_{i,j=1}^nE_{ij}$. Iterating the basic construction horizontally, we obtain the following grid of finite-dimensional algebras
\begin{IEEEeqnarray}{lCl}\label{spin tower}
\begin{matrix}
\Delta_n &\subset & M_n &\subset^{\,e_2} & \Delta_n\otimes M_n &\subset^{\,e_3} & M_n\otimes M_n &\subset & \cdots\cdots\cr
\cup & & \cup & & \cup & & \cup & & \cr
\mathbb{C} &\subset &  u\Delta_n u^* &\subset^{e_2} & \langle u\Delta_nu^*,e_2\rangle &\subset^{e_3} & \langle \langle u\Delta_nu^*,e_2\rangle,e_3\rangle &\subset & \cdots\cdots
\end{matrix}
\end{IEEEeqnarray}
where the notation $\langle\mathcal{A},x\rangle$ for any inclusion of algebras $\mathcal{A\subset B}$ and $x\in\mathcal{B}$ denotes the subalgebra of $\mathcal{B}$ generated by $\mathcal{A}$ and $x$. We obtain a subfactor $R_u\subset R$ with $[R:R_u]=n$, where $R$ is the hyperfinite type $II_1$ factor. This subfactor $R_u\subset R$ is called a {\it spin model} subfactor or {\it Hadamard} subfactor. Not much is known about this class of subfactors, but it is known that they are all irreducible with their second relative commutants abelian. Jones has emphasized their importance on several occasions, see e.g. \cite{Jo2}.
\smallskip

Now, we briefly discuss the vertex model subfactors. A biunitary matrix $w\in M_n\otimes M_k$ induces the following non-degenerate commuting square \cite{JS}
\[
\begin{matrix}
\bbc\otimes M_k &\subset & M_n\otimes M_k \cr \cup &\ & \cup\cr  \mathbb{C} &\subset &  w(M_n\otimes \bbc) w^*\,\,.
\end{matrix}
\]
The horizontal basic construction, like in \Cref{spin tower}, produces a subfactor $R_w\subset R$ with $[R:R_w]=k^2$. This class of subfactors is known as the vertex model subfactors. Note that, unlike the spin model subfactors, this class of subfactors need not be irreducible.

\subsection{Outer automorphism and its outer period}\label{Sec2.4}

The final ingredient needed in this article is the outer automorphisms and their invariants. Denote by $\text{Aut}(R)$ the group of all automorphisms of the hyperfinite $II_1$ factor $R$. Let Inn$(R)$ be the normal subgroup of inner automorphisms of $R$. An automorphism $\alpha\in\mathrm{Aut}(R)$ is said to be an outer automorphism if $\alpha\notin\mathrm{Inn}(R)$. Two automorphisms $\alpha,\beta \in \text{Aut}(R)$ are said to be outerly conjugate if $\sigma \alpha \sigma^{-1}$ and $\beta$ belong to the same equivalence class of $\text{Out}(R):=\text{Aut}(R)/\text{Inn}(R)$.

In \cite{Connes}, the outer invariants are defined for automorphisms that characterize them up to outer conjugacy. For a given automorphism $\alpha \in \text{Aut}(R)$, if there exist $p_{0}(\alpha) \in \bbn$ satisfying the following
\begin{equation}\label{outer period}
\{ n \in \mathbb{Z} \, : \, \alpha^{n} \in \text{Inn}(R)\}=p_{0}(\alpha) \mathbb{Z},
\end{equation}
\noindent then $p_{0}(\alpha)$ is called the outer period of  $\alpha$. Note that $p_{0}(\alpha)$ is the smallest element of $\bbn$ satisfying $\alpha^{p_{0}(\alpha)}=\text{Ad}_{u}$ for some unitary $u\in R$. Further note that there exists $\gamma \in \bbc$ of modulus $1$ satisfying $\alpha(u)=\gamma u$. The pair $(p_{0}(\alpha), \gamma)$ is the outer invariant of $\alpha$. If no such $p_{0}(\alpha)$ exists that satisfies \Cref{outer period}, we say the outer period of $\alpha$ is zero.


\newsection{Construction of a family of biunitary matrices}\label{Sec3}

Let $n$ be a natural number and consider the direct system $\{(\mathscr{A}_j,\phi_j):j\in\bbn\}$, where $\mathscr{A}_j:=M_n^{\otimes\,j}\cong M_{n^j}$, and for each $j\in\bbn,\,\phi_j:\mathscr{A}_j\hookrightarrow\mathscr{A}_{j+1}$ is given by the unital algebra embedding $x\mapsto 1_n\otimes x$ with $1_n$ being the identity matrix in $M_n$. Then, $\cup_j\mathscr{A}_j$ is naturally a unital algebra, and by abuse of notation, we denote it by $M_\infty$. Thus, in our context, $M_\infty:=\cup_{j=1}^\infty\,M_{n^j}$, where we omit the injective maps $\phi_j$ for obvious reason. An element $u\in M_\infty$ shall be called `unitary' if $u$ is a unitary matrix in some $M_{n^j}$. We denote $\mathcal{U}(M_\infty):=\cup_{j=1}^\infty\,U(n^j)$, where $U(n^j)$ is the group of unitary matrices of order $n^j$. Similarly, we denote by $\mathcal{BU}(M_\infty)$, the set of all biunitary matrices in $M_\infty$ (recall from \Cref{Sec2}). By definition, $\mathcal{BU}(M_\infty)\subseteq\mathcal{U}(M_\infty)$ and the inclusion is strict. Recall that $\mathbb{H}(M_n)$ denotes the set of all complex Hadamard matrices in $M_n$. For any element $u=(u_{ij})\in\mathbb{H}(M_n)$, consider the following element
\[
D_u = \sqrt{n} \sum_{i=1}^n \sum_{j=1}^n \overline{u}_{ij} (E_{ii} \otimes E_{jj})
\]
in $U(n^2)$ associated to $u$. Here, $E_{ij}$ denotes the matrix units in $M_n$. With the convention $u_0:=u$, it is easy to verify that the following recursive formula
\begin{eqnarray}\label{tower unitary}
u_{2k+1} &:=& (I_n \otimes u_{2k})(D_u \otimes I_n^{(k)}),\nonumber\\
u_{2k} &:=& u_{2k-1}(u \otimes I_n^{(k)}),
\end{eqnarray}
where $k$ varies over $\bbn\cup\{0\}$, induces a sequence of elements in $\mathcal{U}(M_\infty)$. Let $\prod_{j=1}^\infty\mathcal{U}(M_\infty)$ denote the product in the category of sets. Thus, any arbitrary element of $\prod_{j=1}^\infty\mathcal{U}(M_\infty)$ is a sequence in $\mathcal{U}(M_\infty)$. Therefore, for each $n\in\bbn$, \Cref{tower unitary} determines the following map
\begin{center}
$\Phi_n:\mathbb{H}(M_n)\longrightarrow\prod_{j=1}^\infty\mathcal{U}(M_\infty)$
\end{center}
given by
\begin{eqnarray}\label{the map}
u\mapsto\Phi_n(u):=u_1\oplus u_2\oplus\ldots
\end{eqnarray}
where the unitary matrices $u_j$'s are determined recursively starting with $u$.

\begin{rmrk}\rm
The unitary matrices $u_j$'s appear naturally in the context of Hadamard subfactors, see \cite{N} for instance. Note that, in particular, $u_{2k}\in M^{(k+1)}_n$ and $u_{2k+1}\in\Delta_n\otimes M^{(k+1)}_n\subseteq M^{(k+2)}_n$.
\end{rmrk}

Fix $n\in\bbn$ and take $u,w\in\mathbb{H}(M_n)$. We have the corresponding elements $\Phi_n(u)$ and $\Phi_n(w)$ in $\prod_{j=1}^\infty\mathcal{U}(M_\infty)$ obtained in \Cref{the map}. In other words, we have the following map
\begin{align*}
(\Phi_n,\Phi_n):\mathbb{H}(M_n)\times\mathbb{H}(M_n)\longrightarrow&\prod_{j=1}^\infty\mathcal{U}(M_\infty)\times\prod_{j=1}^\infty\mathcal{U}(M_\infty)\nonumber\\
(u,w)\longmapsto&(\Phi_n(u),\Phi_n(w))\,.
\end{align*}
Recall the set $\mathcal{BU}(M_\infty)$, the set of all biunitary matrices in $M_\infty$. We denote the $j$-fold tensor product $M_n^{\otimes j}$ of $M_n$'s by $M_n^{(j)}$ for notational clarity throughout the article. For each $\ell\in \bbn \cup \{0\}$, consider the unitary matrices $V_\ell:=\sum_{i,j}E_{ij} \ot E_{ji} \ot I^{(\ell)}_n$ in $M^{(\ell+2)}_n$, that is, $V_\ell\in\mathcal{U}(M_\infty)$ for each $\ell$, where $E_{ij}$ denotes the matrix units. Note that each $V_\ell$ flips the first two legs in $M_n\otimes M_n\otimes M_n^{(\ell)}$, and acts as identity on $M_n^{(\ell)}$.
\smallskip

The following theorem is the main goal of this section.

\begin{thm}\label{theorem1}
Fix $n\in\bbn$ and let $u,w\in\mathbb{H}(M_n)$ with $u\neq w$. Consider the unitary matrices $V_\ell:=\sum_{i,j}E_{ij} \ot E_{ji} \ot I^{(\ell)}_n\in M^{(\ell+2)}_n$. Then, the unitary matrices $\,BU(u,w;\ell):=u_{2\ell+2}w_{2\ell+1}V_\ell$ in $M_n^{(\ell+2)}$ defined in \Cref{tower unitary} are biunitary matrices for each $\ell\in\bbn\cup\{0\}$.
\end{thm}

In other words, we obtain the following map
\begin{align*}
\Phi:\mathbb{H}(M_n)\times\mathbb{H}(M_n)&\longrightarrow\prod_{j=1}^\infty\mathcal{BU}(M_\infty)\subseteq\prod_{j=1}^\infty\mathcal{U}(M_\infty)\\
(u,w)&\longmapsto BU(u,w;0)\oplus BU(u,w;1)\oplus\ldots\,.
\end{align*}
This says that starting with two fixed complex Hadamard matrices of the same order, one can obtain a sequence of biunitary matrices of distinct orders, and the order of each element of the sequence increases to infinity.
\smallskip

For transparency, we break the proof into a few lemmas. First, we fix some notations.
\begin{notation}\label{not}\rm
\begin{enumerate}[$(i)$]
\item For $k \in \mathbb{N} \cup \{0\}$ and $u\in\mathbb{H}(M_n)$, we let
\begin{align*}
A_{2k}:=M_n\otimes M^{(k)}_n\,,\quad& A_{2k+1}:=\Delta_n\otimes M_n \otimes M^{(k)}_n\,,\\
B^{u}_{2k}:=\mbox{Ad}_{u_{2k}}\big(\Delta_n\otimes M_n^{(k)}\big)\,,\quad& B^{u}_{2k+1}:=\mbox{Ad}_{u_{2k+1}}\big(M_n \otimes M_n^{(k)}\big)\,,
\end{align*}
with respect to the natural embedding $M^{(k)}_n\hookrightarrow\Delta_n\otimes M^{(k)}_n$ from the left, where `$\mbox{Ad}$' denotes the adjoint representation.
\item We fix the following projection matrices
\[
e_1:=\sum_{i,j=1}^nE_{ij}\quad,\quad e_2:=\sum_{i=1}^nE_{ii}\otimes E_{ii}
\]
in $M_n$ and $M_{n^2}$ respectively. We also define the following sequence of projections
\[
e_{2k+1}=e_1\otimes I^{(k)}_n\quad,\quad e_{2k+2}=e_2 \otimes I^{(k)}_n
\] for $k\in\bbn$. Note that $\{e_1,e_2,e_{2k+1},e_{2k+2}:k\in\bbn\}\subseteq M_\infty$.
\end{enumerate}
\end{notation}

\begin{lmma}\label{unitary projection relation}
Fix $n\in\bbn$ and let $u\in\mathbb{H}(M_n)$. Consider the element $\Phi_n(u)$ in $\prod_{j=1}^\infty\mathcal{U}(M_\infty)$ obtained in \Cref{the map}. Then, $u_ke_ku_k^*=e_{k+1}$ for all $k\in\bbn$, where the projections $e_k$'s are as defined in \Cref{not}.
\end{lmma}
\begin{prf}
For $k=1$, we have
\begin{IEEEeqnarray*}{lCl}
\mathrm{Ad}_{u_1}(I_n \otimes e_1) &=& \frac{1}{n} \text{Ad}_{(I_n \otimes u)} \text{Ad}_{D_u} \Big( I_n \otimes \sum_{i,j} E_{i,j} \Big)\\
&=& \mathrm{Ad}_{(I_n \otimes u)} \Big( \sum_{i,j,k} \overline{u_{k,i}} u_{j,k} \left( E_{k,k} \otimes E_{i,j} \right) \Big)\\
&=& \mathrm{Ad}_{(I_n \otimes u)} \left( \text{Ad}_{(I_n \otimes u^*)}(e_2) \right)\\
&=& e_2.
\end{IEEEeqnarray*}
The general case follows by induction on $k$, and we omit the details.\qed
\end{prf}

\begin{rmrk}\rm
The above lemma first appeared in \cite{N} in a different disguise.
\end{rmrk}

For each $k \in \bbn \cup \{0\}$, consider the following subalgebras of $M_\infty:$
\begin{IEEEeqnarray}{lCl}\label{notation two}
N^{u,w}_{2k} &:=& \text{Ad}_{u_{2k+1}(I_n \otimes w_{2k})}(I_n \otimes \dt \otimes I_n^{(k)}),\nonumber\\
N^{u,w}_{2k+1} &:=& \text{Ad}_{u_{2k+2}w_{2k+1}}(I_n \otimes M_n \otimes I_n^{(k)}).
\end{IEEEeqnarray}
\begin{lmma}
For each $k \in \bbn \cup \{0\}$, the following quadruple
\begin{equation}
\mathscr{I}_k:=
\begin{matrix}\label{infinite quadruples}
\mathbb{C} \otimes M_n^{(k+1)} & \subset & \dt \otimes M_n^{(k+1)} \\
\cup & & \cup \\
\mathbb{C} & \subset & N_{2k}
\end{matrix}
\end{equation}
is a non-degenerate commuting square.
\end{lmma}
\begin{prf}
From \Cref{tower unitary}, we know that the unitary matrix $w_{2k} \in M_n \otimes M_n^{(k)}$, and hence,
\[
N^{u,w}_{2k} = \text{Ad}_{u_{2k+1}(I_n \otimes w_{2k})}(I_n \otimes \dt \otimes I_n^{(k)}) \subset B_{2k+1}^u.
\]
Therefore, we have
\begin{equation}\label{conditional}
E_{N^{u,w}_{2k}}(\bbc \otimes M_n^{(k+1)}) = E_{N^{u,w}_{2k}} E_{B_{2k+1}^u}(\bbc \ot M^{(k+1)}_n).
\end{equation}
Since $u$ is a complex Hadamard matrix, the following quadruple
\begin{equation*}
\begin{matrix}
B^{u}_{2k+1} & \subset & \dt \otimes M_n^{(k+1)} \\
\cup & & \cup \\
B^{u}_{2k} & \subset &  \bbc \ot M^{(k+1)}_{n}
\end{matrix}
\end{equation*}
is a non-degenerate commuting square. Consequently, $E_{B_{2k+1}^u}(\bbc \otimes M_n^{(k+1)}) = B_{2k}^u$. Substituting in \Cref{conditional}, we obtain
\begin{equation}\label{commute exp}
 E_{N^{u,w}_{2k}} (\bbc \otimes M_n^{(k+1)}) = E_{N^{u,w}_{2k}}(B_{2k}^u)   
\end{equation}
From the definition of $B_{2k}^u$, we have
\begin{equation*}
B_{2k}^u = \text{Ad}_{(I_n \otimes u_{2k})}(\bbc \otimes \dt \otimes M_n^{(k)}),
\end{equation*}
Since $u_{2k+1} = (I_n \otimes u_{2k})(D_u \otimes I_n^{(k)})$ from \Cref{tower unitary}, we can rewrite the above as
\begin{equation}\label{k-th space}
B_{2k}^u = \text{Ad}_{u_{2k+1}}(\bbc \otimes \dt \otimes M_n^{(k)}).
\end{equation}
Therefore, substituting the definition of $N^{u,w}_{2k}$ from \Cref{notation two}, we get
\begin{eqnarray}\label{equation of n2k}
E_{N^{u,w}_{2k}}(B_{2k}^u) &= &\text{Ad}_{u_{2k+1}} \big\{E_{\mathscr{A}}(B^{2k}_u)\big\}  \hspace{2mm}(\text{by \Cref{k-th space}}) \nonumber\\
&=&\text{Ad}_{u_{2k+1}}\big(E_{\mathscr{A}}(I_n \otimes \dt \otimes M_n^{(k)})\big),
\end{eqnarray}
where we denote $\mathscr{A}:=\text{Ad}_{(I_n \otimes w_{2k})}(I_n \otimes \dt \otimes I_n^{(k)})$.
\smallskip

\noindent \textbf{Claim}:~ $E_{\mathscr{A}}(I_n \otimes \dt \otimes M_n^{(k)})=\bbc$.
\smallskip

\noindent To establish the claim, first observe that from \Cref{tower unitary}, we have $w_{2k} = w_{2k-1}(w \otimes I_n^{(k)})$, where $w_{2k-1} \in \dt \otimes M_n^{(k)}$. Therefore, it follows that
\begin{equation}\label{mid equation}
 \text{Ad}_{(I_n \otimes w_{2k}^*)}(I_n \otimes \dt \otimes M_n^{(k)})=(I_n \ot w^{*} \dt w \ot M^{(k)}_n)   
\end{equation}
where $w^{*}$ is the adjoint of $w$. Since $\mathscr{A} = \text{Ad}_{(I_n \otimes w_{2k})}(I_n \otimes \dt \otimes I_n^{(k)})$, using \Cref{mid equation} we establish our claim.

So from \Cref{equation of n2k}, it follows that $E_{N^{u,w}_{2k}}(B_{2k}^{u})=\bbc$, and therefore, by \Cref{commute exp} we conclude that  $E_{N^{u,w}_{2k}} (\bbc \otimes M_n^{(k+1)})=\bbc$. Thus the quadruple $\mathscr{I}_k$ is a commuting square. Since the norms of the inclusion matrices for both the vertical embeddings are the same, we conclude that $\mathscr{I}_k$ is a non-degenerate commuting square, which completes the proof.\qed
\end{prf}

\begin{lmma}\label{subset of algebra}
For $k\in\bbn\cup\{0\}$, consider the subalgebra $\langle N^{u,w}_{2k}, e_{2k+3}\rangle$ of $M_\infty$, generated by $N^{u,w}_{2k}$ and the projection $e_{2k+3}$ defined in \Cref{not}. Then, we have the following:
\[
N^{u,w}_{2k+1}=\langle N^{u,w}_{2k}, e_{2k+3} \rangle\qquad\mbox{and}\qquad N^{u,w}_{2k+2} \subseteq\langle N^{u,w}_{2k+1}, \,e_{2k+4}\rangle\,.
\]
\end{lmma}
\begin{prf}
From \Cref{tower unitary}, we have the following identities:
\[
u_{2k+2}=u_{2k+1} (u \ot I^{(k+1)}_n) \quad \text{and} \quad w_{2k+1}=(I_n \ot w_{2k})(D_w \ot I^{(k)}_n)
\]
where $D_w\in\Delta_n\ot\Delta_n$. Using this fact, we see the following
\[\text{Ad}_{u_{2k+2}}\text{Ad}_{(I_n \ot w_{2k})}(I_n \ot \dt \ot I^{(k)}_n)=\text{Ad}_{u_{2k+1}}(\text{Ad}_{(I_n \ot w_{2k})}(I_n \ot \dt \ot I^{(k)}_n))
\]
\[
\text{Ad}_{(I_n \ot w_{2k})}(I_n \ot \dt \ot I^{(k)}_n)=\text{Ad}_{w_{2k+1}}(I_n \ot \dt \ot I^{(k)}_n)\,.
\]
Therefore, it follows that
\begin{eqnarray}\label{rewrite one}
\text{Ad}_{u_{2k+2}w_{2k+1}}(I_n \ot \dt \ot I^{(k)}_n)&=&\text{Ad}_{u_{2k+1}}(\text{Ad}_{(I_n \ot w_{2k})}(I_n \ot \dt \ot I^{(k)}_n))=N_{2k}\,.
\end{eqnarray}
Along the same line, again by using \Cref{tower unitary}, one can prove the following:
\begin{eqnarray}\label{rewrite two}
\text{Ad}_{u_{2k+2}w_{2k+1}}(I_n \ot M_n \ot I^{(k)}_n)&=&\text{Ad}_{u_{2k+3}(I_n \ot w_{2k+2})}(I_n \ot M_n \ot I^{(k)}_n)=N_{2k+1}\,.
\end{eqnarray}
Now,
\begin{eqnarray}\label{prove one}
& & \langle N^{u,w}_{2k}, \,e_{2k+3}\rangle\nonumber\\
&=&N^{u,w}_{2k}e_{2k+3}N^{u,w}_{2k} \nonumber \hspace{2mm} (\text{by \, \Cref{rewrite one}})\\
&=&\text{Ad}_{u_{2k+2}w_{2k+1}}(I_n \ot \dt \ot I^{(k)}_n) e_{2k+3} \text{Ad}_{u_{2k+2}w_{2k+1}}(I_n \ot \dt \ot I^{(k)}_n)\nonumber\\
&=&\text{Ad}_{u_{2k+2}w_{2k+1}}\{(I_n \ot \dt \ot I^{(k)}_n)  \{\text{Ad}_{(u_{2k+2}w_{2k+1})^{*}}e_{2k+3}\}(I_n \ot \dt \ot I^{(k)}_n) \}\,.
\end{eqnarray}
Using \Cref{unitary projection relation}, one obtains $\text{Ad}_{(u_{2k+2}w_{2k+1})^{*}}e_{2k+3}=e_{2k+1}$, and therefore \Cref{prove one} becomes the following:
\begin{IEEEeqnarray*}{lCl}
& & \langle N^{u,w}_{2k}, \,e_{2k+3}\rangle\\
&=&\text{Ad}_{u_{2k+2}w_{2k+1}}\{(I_n \ot \dt \ot I^{(k)}_n) e_{2k+1}(I_n \ot \dt \ot I^{(k)}_n) \}\\
&=&\text{Ad}_{u_{2k+2}w_{2k+1}}\{(I_n \ot \dt \ot I^{(k)}_n)(I_n \ot e_1 \ot I^{(k)}_n)(I_n \ot \dt \ot I^{(k)}_n)\}\\
&=&\text{Ad}_{u_{2k+2}w_{2k+1}}(I_n \ot M_n \ot I^{(k)}_n) \hspace{4mm}(\text{by \Cref{notation two}})\\
&=&N^{u,w}_{2k+1}\,.
\end{IEEEeqnarray*}
This completes the first part of the proof.

Finally, along the similar lines, using \Cref{unitary projection relation} and \Cref{rewrite two}, one can prove the following\,:
\[
\langle N^{u,w}_{2k+1}, e_{2k+4}\rangle=\text{Ad}_{u_{2k+3} (I_n \ot w_{2k+2})}(I_n \ot \dt \ot M_n \ot I^{(k)}_n)\,,
\]
which proves that $N^{u,w}_{2k+3}\subseteq\langle N^{u,w}_{2k+1}, \,e_{2k+4}\rangle$.\qed
\end{prf}
\medskip

\noindent\textbf{Proof of the \Cref{theorem1}:}
First, we observe the following
\begin{eqnarray}\label{biunitary}
   \langle N^{u,w}_{2k}, e_{2k+3} \rangle &=&N^{u,w}_{2k+1}   \hspace{2mm} (\text{by \,\Cref{subset of algebra}})\nonumber \\
    &=& \text{Ad}_{u_{2k+2}w_{2k+1}}(I_n \otimes M_n \otimes I_n^{(k)}) \hspace{2mm} (\text{by \Cref{notation two}})\nonumber\\
    &=& \text{Ad}_{u_{2k+2}w_{2k+1}V_k}(M_n \otimes I_n \otimes I_n^{(k)})
\end{eqnarray}
where $V_k:=\sum_{i,j}E_{ij}\ot E_{ji} \ot I^{(k)}_n$ is a unitary matrix. The following quadruples
\begin{IEEEeqnarray*}{lCl}
\begin{matrix}
\bbc \ot M^{(k+1)}_{n} &\subset & \dt \otimes M^{(k+1)}_{n} & \subset^{\,e_{2k+3}} &  M_n \ot M^{(k+1)}_n\cr \cup &\ &\cup  & & \cup \cr \bbc &\subset & N^{u,w}_{2k} & \subset^{\,e_{2k+3}} &  \langle N^{u,w}_{2k}, \, e_{2k+3} \rangle
\end{matrix}
\end{IEEEeqnarray*}
is the basic construction of the non-degenerate commuting square $\mathscr{I}_k$ in \Cref{infinite quadruples} with the Jone's projection $e_{2k+3}=e_1 \ot I^{(k+1)}_n$. Since the following quadruple
\begin{eqnarray*}
\mathscr{L}_k&=& \begin{matrix}
\bbc \ot M^{(k+1)}_{n} &\subset & M_n \otimes M^{(k+1)}_{n} \cr \cup &\ &\cup\cr \bbc &\subset & \langle N^{u,w}_{2k}, e_{2k+3} \rangle
\end{matrix} 
\end{eqnarray*}
is a non-degenerate commuting square, the unitary matrix $u_{2k+2}w_{2k+1}V_k$ in \Cref{biunitary} must be a biunitary matrix by \Cref{biunitary thm}. This concludes the proof.\qed


\newsection{Application to Subfactors}\label{Sec4}

Traversing from finite to infinite dimension, in this section we apply our investigations in \Cref{Sec3} to subfactor theory. Starting with an element $(u,w)$ in $\mathbb{H}(M_n)\times \mathbb{H}(M_n)$, we construct a nested sequence of vertex model subfactors with distinct Jones indices in the hyperfinite type $II_1$ factor $R$ such that the Hadamard/spin model subfactor $R_u\subset R$ is an intermediate subfactor to each of them. Moreover, this sequence of subfactors is not the tower of downward basic construction. Our construction is noncommutative; that is, subfactors arising from $(u,w)$ need not be isomorphic with that arising from $(w,u)$.

\subsection{A family of vertex model subfactors}\label{Sec4.1}

Let $u,w$ be two complex Hadamard matrices of the same order $n$. We have a sequence of biunitary matrices $\{BU(u,w;k)\}_{k=0}^\infty$ obtained to \Cref{theorem1}. Using itereated basic construction, we construct a nested sequence of vertex model subfactors from $\{BU(u,w;k)\}_{k=0}^\infty$. More precisely, we prove the following theorem.

\begin{thm}\label{theorem2}
Fix any $n\geq 2$. Let $u\in\mathbb{H}(M_n)$ be a complex Hadamard matrix and $R_u \subset R$ be the corresponding Hadamard/spin model subfactor. Choose any $w\in\mathbb{H}(M_n)$. For each $k\in\mathbb{N}\cup\{0\}$, there exists a vertex model subfactor $R^{u, w}_k \subset R$ such that
\begin{enumerate}[$(i)$]
\item $R^{u, w}_k\subset R_u$, and the Jones index is given by $[R:R^{u, w}_k]=n^{2(k+1)}$;
\item The subfactors $\{R^{u, w}_k\subset R:k\in\mathbb{N}\cup\{0\}\}$ satisfy the inclusion $R^{u, w}_{k+1} \subset R^{u, w}_k$ for each $k$.
\end{enumerate}
Therefore, we have a nested sequence of vertex model subfactors
\[
\cdots \subset R^{u, w}_{k+1}\subset R^{u, w}_k\subset\cdots\subset R^{u, w}_0\subset R
\]
each having the Hadamard/spin model subfactor $R_u\subset R$ as an intermediate subfactor.
\end{thm}
\begin{prf}
We divide the proof into three steps.

\noindent \textbf{Step 1:} $R^{u,w}_k\subset R$ is a vertex model subfactor with $[R : R_{k}^{u,w}] = n^{2(k+1)}$.

\noindent For each $k \in \mathbb{N} \cup \{0\}$, we have the following quadruple
\begin{eqnarray*}
\mathscr{I}_{k}&=& \quad \begin{matrix}
\bbc \ot M^{(k+1)}_{n} &\subset & \dt \otimes M^{(k+1)}_{n} \cr \cup &\ &\cup\cr \bbc &\subset{} & N^{u,w}_{2k}
\end{matrix} 
\end{eqnarray*}
which is a non-degenerate commuting square. Using the iterated basic construction to the above non-degenerate commuting square $\mathscr{I}_k$, we obtain the following grid of finite-dimensional algebras\,:
\begin{eqnarray}\label{construction middle vertex model}
\begin{matrix}
\bbc \ot M^{(k+1)}_n  &\subset & \Delta_n\otimes M^{(k+1)}_n &\subset^{\,e_{2k+3}}  & \subset M_n\otimes M^{(k+1)}_n &\subset^{e_{2k+4}}  \cdots &\subset R\cr
\cup & & \cup & & \cup &  & \,\, \, \cup \cr
\bbc  &\subset & N^{u,w}_{2k} &\subset^{\,e_{2k+3}}  & \subset \langle N^{u,w}_{2k} \, ,e_{2k+3}\rangle & \subset^{e_{2k+4}}  \cdots \cdots & \, \subset R_{k}^{u,w}.
\end{matrix}
\end{eqnarray}
Here, the projections $e_{2k+2i+1} = e_1 \otimes I_n^{(k+i)}$ and $e_{2k+2i+2} = e_2 \otimes I_n^{(k+i)}$, for all $i \in \mathbb{N}$, are the Jones' projections. Consider the following von Neumann algebra  
\begin{eqnarray}\label{middle vertex subfactor}
R_{k}^{u,w}:=\overline{\bigvee \{N^{u,w}_{2k}, \,e_{2k+2+i} \,: \, i\in \bbn\}}^{\,\text{sot}}. 
\end{eqnarray}
It follows that $R_{k}^{u,w}$ is a $II_1$-subfactor of the hyperfinite factor $R$. Finally, the inclusion matrix for the inclusion $\bbc \subset M_n^{(k+1)}$ is the $1\times 1$ matrix $[n^{(k+1)}]$, and hence by \Cref{prop}, $[R : R_{k}^{u,w}] = n^{2(k+1)}$. Now, from \Cref{construction middle vertex model} the following non-degenerate commuting square
\begin{eqnarray*}
\mathscr{L}_{k}&=& \quad \begin{matrix}
\bbc \ot M^{(k+1)}_{n} &\subset & M_n \otimes M^{(k+1)}_{n} \cr \cup &\ &\cup\cr \bbc &\subset & \langle N^{u,w}_{2k}, e_{2k+3} \rangle
\end{matrix} 
\end{eqnarray*}
is a basic construction of the nondegenerate commuting square $\mathscr{I}_k$. Here, $\langle N^{u,w}_{2k}, e_{2k+3}\rangle =N^{u,w}_{2k+1}=\mathrm{Ad}_{u_{2k+2}w_{2k+1}V_k}(M_n \otimes I_n \otimes I_n^{(k)})$ (see \Cref{biunitary}), and $BU(u,w;k):=u_{2k+2}w_{2k+1}V_k$ is a biunitary matrix of order $n^{k+1}$ (see \Cref{theorem1}). Thus, $R^{u,w}_k\subset R$ is a vertex model subfactor. This completes Step $1$.
\smallskip

\noindent \textbf{Step 2:} $R_{k+1}^{u,w}\subset R_{k}^{u,w}$ for all $k\in\mathbb{N} \cup \{0\}.$

\noindent Note that by \Cref{middle vertex subfactor}, we have
\begin{eqnarray}\label{next vertex}
R_{k+1}^{u,w}=\overline{\bigvee \{N^{u,w}_{2k+2}, \,e_{2k+4+i} \,: \, i\in \bbn\}}^{\,\text{sot}}. 
\end{eqnarray}
From \cref{subset of algebra}, we know that $N^{u,w}_{2k+2} \subseteq \langle N_{2k+1}, e_{2k+4} \rangle$, and since $N^{u,w}_{2k+1} = \langle N^{u,w}_{2k}, e_{2k+3} \rangle$, it follows that  
\begin{eqnarray*}
\{N^{u,w}_{2k+2}, \,e_{2k+4+i} : i \in \bbn\}\subseteq \{N^{u,w}_{2k}, \,e_{2k+2+i} : i \in \bbn\}
\end{eqnarray*}
Combining this with \Cref{middle vertex subfactor,next vertex}, we conclude the proof of this step.
\smallskip

\noindent \textbf{Step 3:} $R^{u,w}_{k} \subset R_u\,\,\mbox{for all }k\in\mathbb{N}\cup\{0\}$.

\noindent By Step $2$, it suffices to only show that $R^{u,w}_0 \subset R_u$. From \Cref{middle vertex subfactor}, we observe that
\begin{eqnarray*}
R^{u,w}_0 &=& \overline{\bigvee \{N^{u,w}_0, e_{2+i} : i \in \mathbb{N} \}}^{\,\text{sot}},
\end{eqnarray*}
where $N^{u,w}_0 = \text{Ad}_{u_1}(\mathbb{C} \otimes w \Delta_n w^*)$. Additionally, from \Cref{not}, recall that 
$R_u = \overline{\bigcup_{j} B^u_j}^{\text{sot}}$, where
$B^u_1 = \text{Ad}_{u_1}(\mathbb{C} \otimes M_n)$. This implies that 
$N^{u,w}_0 \subset B^u_1 \subset R_u$. Furthermore, we have
$e_i \in R_u$ for all $i \geq 2$. Consequently, we deduce that
\begin{eqnarray*}
 \bigvee \{N^{u,w}_0, e_{2+i} : i \in \mathbb{N} \} \subset R_u.
\end{eqnarray*}
This establishes the fact that $R^{u,w}_0\subset R_u$.\qed
\end{prf}

At this stage, the following natural questions arise.
\smallskip

\noindent\textbf{Question A:} Is the vertex model subfactor $R^{u, w}_0\subset R$ the downward basic construction of the Hadamard subfactor $R_u\subset R$?
\smallskip

\noindent\textbf{Question B:} Is the nested sequence of vertex model subfactors $\,\,\cdots\cdots\subset R^{u, w}_1\subset R^{u, w}_0\subset R$ a tower of downward basic construction of $R^{u, w}_0\subset R$?
\smallskip

If the answer to any of the above questions is always yes, then our construction would not be considered as new. Below, we show that this is not the case.
\smallskip

\noindent\textbf{Answer to Question A:} Let $n\geq 4$ be an even number and $\omega$ be a primitive $n$-th root of unity. Consider a permutation matrix $P\in M_n$ such that $\mathrm{Ad}_{PF_n}(\mathscr{D}_{\frac{n}{2}})\neq\alpha\,\mathrm{Ad}_{\mathscr{D}_1PF_n}(\mathscr{D}_{\frac{n}{2}})$ for any $\alpha\in\bbc$, where
\begin{eqnarray*}
\mathscr{D}_1&:=& \mathrm{diag}\{1, \omega, \omega^2, \cdots, \omega^{n-1}\},\nonumber\\
\mathscr{D}_k &:=& \mathscr{D}_1^k=\mathrm{diag}\{1, \omega^k, \omega^{2k}, \cdots, \omega^{(n-1)k}\}.
\end{eqnarray*}
Such a $P$ always exists, for example, choose $P:=E_{13}+E_{21}+E_{32}+\sum_{k=4}^{n} E_{kk}$, where $E_{ij}$'s are the matrix units. Take $u=F_n$ and $w=PF_n$. We show that the vertex model subfactor $R^{F_n,PF_n}_0 \subset R$ is irreducible, which immediately says that $R^{F_n, PF_n}_0\subset R_{F_n}$ can not be the downward basic construction of the spin model subfactor $R_{F_n}\subset R$.

From \Cref{construction middle vertex model}, we have the following non-degenerate commuting square (putting $k=0$)
\[
\begin{matrix}
\bbc \ot M_{n} &\subset & \dt \otimes M_{n} \cr \cup &\ &\cup\cr \bbc &\subset{} & N^{F_n,PF_n}_{0}
\end{matrix} 
\]
where $N^{F_n,PF_n}_0=\text{Ad}_{(F_n)_1(I_n \otimes PF_n)}(\bbc\otimes\Delta_n)$. By construction, the vertex model subfactor $R^{F_n,PF_n}_0 \subset R$ is obtained from the above non-degenerate commuting square using the iterated basic construction. Hence, using the fact that $(F_n)_1=\mbox{bl-diag}\{F_n,F_n\mathscr{D}_{n-1},\ldots,F_n\mathscr{D}_1\}$, by the Ocneanu compactness we have the following\,:
\begin{IEEEeqnarray}{lCl}\label{changed}
(R^{F_n,PF_n}_0)^\prime\cap R &=& (N^{F_n,PF_n}_0)^{'}\cap (\bbc \ot M_n)\nonumber\\
&=& \big(\mathrm{Ad}_{(F_n)_1(I_n \ot PF_n)}(\bbc \ot \Delta_n)\big)^\prime\cap (\bbc \ot M_n)\nonumber\\
&=& \mathrm{Ad}_{F_n}\Big(\bigcap_{k=1}^{n}\mathrm{Ad}_{\mathscr{D}_kPF_n}(\Delta_n)\Big).
\end{IEEEeqnarray}
We claim that $\bigcap_{k=1}^n\mathrm{Ad}_{\mathscr{D}_{k} PF_n}(\Delta_n)=\bbc$. On contrary, assume that $\cap_{k=1}^n\mathrm{Ad}_{\mathscr{D}_kPF_n} (\Delta_n)\neq\bbc$. Then, there exist non-trivial diagonal matrices $D_0, D_1$ such that  $\mathrm{Ad}_{PF_n} (D_0)=\mathrm{Ad}_{\mathscr{D}_1PF_n} (D_1)$. Consider $D_0=\sum_{i=0}^{n-1} r_i \mathscr{D}_i$ and $D_1=\sum_{i=0}^{n-1} s_i \mathscr{D}_i$, and deduce the following  
\begin{equation}\label{last eqn}
r_i \mathrm{Ad}_{PF_n}(\mathscr{D}_i)=s_i\mathrm{Ad}_{\mathscr{D}_1PF_n}(\mathscr{D}_i)\quad\text{for all }\,i.
\end{equation}
Since $D_0\notin\bbc I_n$, there exists some $i_0\in\{1,\ldots,n-1\}$ such that $r_{i_0}\neq 0\,\text{and}\,s_{i_0}\neq 0$. Therefore, it follows from \Cref{last eqn} that $\mathrm{Ad}_{PF_n}(\mathscr{D}_{i_0})=\frac{s_{i_0}}{r_{i_0}}\mathrm{Ad}_{\mathscr{D}_1PF_n}(\mathscr{D}_{i_0})$. Since there exists some $k\in\bbn$ such that $\mathscr{D}_{i_0}^k=\mathscr{D}_{\frac{n}{2}}$, and $n$ is even, we obtain the following\,:
\begin{IEEEeqnarray*}{lCl}
\mathrm{Ad}_{PF_n}(\mathscr{D}_{\frac{n}{2}}) &=& (\mathrm{Ad}_{PF_n}(\mathscr{D}_{i_0}))^{k}\nonumber\\
&=& \Big(\frac{s_{i_0}}{r_{i_0}}\Big)^k(\mathrm{Ad}_{\mathscr{D}_1PF_n}(\mathscr{D}_{i_0}))^{k}\nonumber\\
&=& \Big(\frac{s_{i_0}}{r_{i_0}}\Big)^k\mathrm{Ad}_{\mathscr{D}_1PF_n}(\mathscr{D}_{\frac{n}{2}}).
\end{IEEEeqnarray*}
This contradicts the choice of $P$ \big(with $\alpha=(s_{i_0}/r_{i_0})^k$ in equation $\mathrm{Ad}_{PF_n}(\mathscr{D}_{\frac{n}{2}})\neq\alpha\,\mathrm{Ad}_{\mathscr{D}_1PF_n}(\mathscr{D}_{\frac{n}{2}})$\big). Hence, from \Cref{changed} we conclude that
\begin{equation}\label{to be used next}
 (R^{F_n,PF_n}_0)^{'} \cap R=\bbc.
\end{equation}

\noindent\textbf{Answer to Question B:} Let $F_n$ denote the Fourier matrix of order $n$, and consider $u=w=F_n$. For any $m,k\in\bbn$, the following nested sequence
\[
\cdots\subset R_{(m+1)k-1}^{u,w}\subset R_{mk-1}^{u,w}\subset R_{(m-1)k-1}^{u,w}\subset\cdots\subset R_{2k-1}^{u,w}\subset R_{k-1}^{u,w}\subset R
\]
of vertex model subfactors cannot be obtained as the downward basic construction of the subfactor $R_{k-1}^{u,w}\subset R$. This is because there is no non-zero projection $p\in\big(R^{u,w}_{(m+1)k-1}\big)^\prime\cap R$ satisfying $pxp=E^{R_{mk-1}^{u,w}}_{R_{(m+1)k-1}^{u,w}}(x)p$ for all $x\in R_{mk-1}^{u,w}$. The proof of this fact is provided in the Appendix.


\subsection{Non-commutativity of the construction}\label{Sec4.2}

\begin{figure}
    \centering
    \begin{tikzpicture}
          \draw (-7.8,0) node [scale=1.2] {{$R$}};
          \draw (-8.5, -.8) node  [scale=1.4] [rotate=45]  {$\subset$};
             \draw (-8.5, .8) node [scale=1.4] [rotate=-45] {$\subset$};
          \draw (-9.5, 1.5) node [scale=1.2] {{$R_u$}};
          \draw (-9.5, -1.5) node [scale=1.2] {{$R_w$}};
          \draw (-10.5, -1.5) node [scale=1.4] {{$\subset$}};
         \draw (-10.5, 1.5) node  [scale=1.4] {{$\subset$}};
          \draw (-11.5, 1.5) node [scale=1.2] {{$R^{u,w}_0$}};
           \draw (-11.5, -1.5) node [scale=1.2] {{$R^{w,u}_0$}}; 
            \draw (-12.5, -1.5) node [scale=1.4] {{$\subset$}};
            \draw (-12.5, 1.5) node  [scale=1.4] {{$\subset$}};
            \draw (-13.5, -1.5) node [scale=1.2] {{$R^{w,u}_1$}};
            \draw (-13.5, 1.5) node  [scale=1.2] {{$R^{u,w}_1$}};
         \draw (-14.5, -1.5) node [scale=1.4] {{$\subset$}};
          \draw (-14.5, 1.5) node  [scale=1.4] {{$\subset$}};
           \draw[dashed](-15, -1.5)--(-16, -1.5);
            \draw[dashed](-15, 1.5)--(-16, 1.5); 
            \draw (-16.5, -1.5) node [scale=1.4] {$\subset$};
             \draw (-16.5, 1.5) node [scale=1.4] {$\subset$};
             \draw (-17.5, 1.5) node [scale=1.2] {{$R^{u,w}_k$}};
            \draw (-17.5, -1.5) node [scale=1.2] {{$R^{w,u}_k$}};
            \draw (-18.5, -1.5) node [scale=1.4] {{$\subset$}};
          \draw (-18.5, 1.5) node  [scale=1.4] {{$\subset$}};
           \draw[dashed](-19, -1.5)--(-20, -1.5); 
            \draw[dashed](-19, 1.5)--(-20, 1.5); 
 \end{tikzpicture}
    \caption{Noncommutativity of the construction}\label{fig1}
\end{figure}
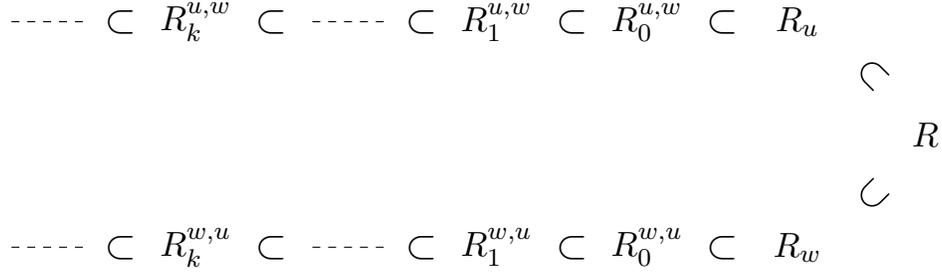

Recall that in the proof of \Cref{theorem2}, we started with $u\in\mathbb{H}(M_n)$ and took help of $w\in\mathbb{H}(M_n)$ to construct the nested sequence of vertex model subfactors
\[
\cdots \subset R^{u, w}_{k+1}\subset R^{u, w}_k\subset\cdots\subset R^{u, w}_0\subset R_u\subset R
\]
that have the Hadamard subfactor $R_u\subset R$ as an intermediate. Instead of this, if we start with $w$ and take the help of $u$ to construct the following nested sequence of vertex model subfactors
\[
\cdots \subset R^{w, u}_{k+1}\subset R^{w, u}_k\subset\cdots\subset R^{w, u}_0\subset R_w\subset R
\]
that have the Hadamard subfactor $R_w\subset R$ as an intermediate, then a natural question arises here whether $R^{u, w}_k\subset R$ are isomorphic to $R^{w, u}_k\subset R$ or not. This situation is depicted in \Cref{fig1}. This is exactly the question of noncommutativity of our construction.

We recall once again the answer to Question A in \Cref{Sec4.1}. Let $n \geq 4$ be an even number and consider the permutation matrix $P$ of order $n$ such that $\mathrm{Ad}_{PF_n}(\mathscr{D}_{\frac{n}{2}})\neq\alpha\,\mathrm{Ad}_{\mathscr{D}_1PF_n}(\mathscr{D}_{\frac{n}{2}})$ for any $\alpha\in\bbc$, where $\mathscr{D}_1=\text{diag}\{1, \omega, \omega^{2},\ldots, \omega^{n-1}\}$ and $\mathscr{D}_{\frac{n}{2}}=\mathscr{D}_{1}^{\frac{n}{2}}$. Take $u=F_n$ and $w=PF_n$. By \Cref{to be used next} we have $(R_0^{u,w})^\prime\cap R=\bbc$. However, it turns out that $(R^{w,u}_0)^\prime\cap R=\bbc^{n}$. This follows as a corollary to a more general result established in \Cref{Sec6} (see \Cref{diagonal subfactor}), and we refrain from exhibiting the proof here. Hence, the vertex model subfactors $R^{u,w}_0 \subset R$ and $R^{w,u}_0 \subset R$ are clearly non-isomorphic, which establishes the noncommutativity of our construction.

As a simple corollary, we obtain that the biunitary matrices $BU(u,w;\ell)$ constructed in \Cref{theorem1} respect the ordered pair $(u,w)$. More precisely, the matrices $BU(u,w;\ell)$ and $BU(w,u;\ell)$ are not unitarily equivalent for any $\ell\in\bbn$, as otherwise, the corresponding vertex model subfactors would have always been isomorphic.


\section{Inner and outer automorphisms}\label{Sec5}

\begin{figure}
\centering
\begin{tikzpicture}
          \draw (-7.8,0) node [scale=1.2] {{$R$}};
          \draw (-8.5, -.8) node  [scale=1.4] [rotate=45]  {$\subset$};
             \draw (-8.5, .8) node [scale=1.4] [rotate=-45] {$\subset$};
          \draw (-9.5,1.5) node [scale=1.2] {{$R_u$}};
          \draw (-9.5,-1.5) node [scale=1.2] {{$R_{\widetilde{u}}$}};
          \draw (-10.5, -1.5) node [scale=1.4] {{$\subset$}};
         \draw (-10.5, 1.5) node  [scale=1.4] {{$\subset$}};
          \draw (-11.5, 1.5) node [scale=1.2] {{$R^{u,v}_0$}};
           \draw (-11.5,-1.5) node [scale=1.2] {{$R^{\widetilde{u},\widetilde{v}}_0$}}; 
            \draw (-12.5, -1.5) node [scale=1.4] {{$\subset$}};
            \draw (-12.5, 1.5) node  [scale=1.4] {{$\subset$}};
            \draw (-13.5, -1.5) node [scale=1.2] {{$R^{\widetilde{u},\widetilde{v}}_1$}};
            \draw (-13.5, 1.5) node  [scale=1.2] {{$R^{u,v}_1$}};
         \draw (-14.5,-1.5) node [scale=1.4] {{$\subset$}};
          \draw (-14.5, 1.5) node  [scale=1.4] {{$\subset$}};
           \draw[dashed](-15, -1.5)--(-16, -1.5); 
            \draw[dashed](-15, 1.5)--(-16, 1.5); 
            \draw (-16.5, -1.5) node [scale=1.4] {$\subset$};
             \draw (-16.5, 1.5) node [scale=1.4] {$\subset$};
             \draw (-17.5, 1.5) node [scale=1.2] {{$R^{u,v}_k$}};
            \draw (-17.5, -1.5) node [scale=1.2] {{$R^{\widetilde{u},\widetilde{v}}_k$}};
            \draw (-18.5, -1.5) node [scale=1.4] {{$\subset$}};
          \draw (-18.5, 1.5) node  [scale=1.4] {{$\subset$}};
           \draw[dashed](-19, -1.5)--(-20, -1.5); 
            \draw[dashed](-19, 1.5)--(-20, 1.5); 
\end{tikzpicture}
\caption{Infinitely many pairs of vertex model subfactors}\label{fig2}
\end{figure}
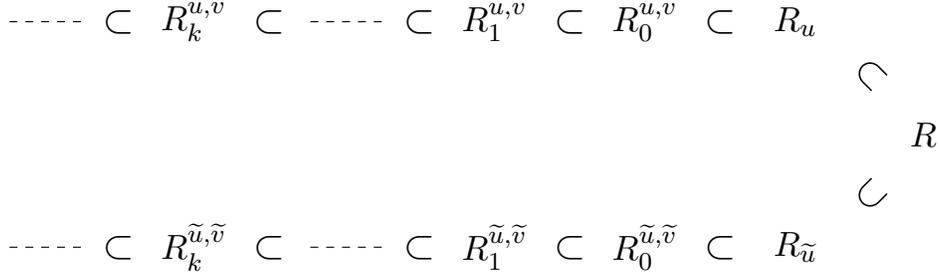

In \Cref{Sec4}, we constructed infinitely many vertex model subfactors using the input data that consists of a pair $(u,v)$ of complex Hadamard matrices of the same order. Consider two such pairs $(u,v),\,(\widetilde{u},\widetilde{v})$. By \Cref{theorem2}, we have the infinitely many pairs of vertex model subfactors $R_k^{u,v}\subset R$ and $R_k^{\widetilde{u},\widetilde{v}}\subset R$ for each $k\in\bbn\cup\{0\}$, depicted in \Cref{fig2}.
The following natural question arises.
\smallskip

\noindent\textbf{Question:} Are the pairs of vertex model subfactors $R_k^{u,v}\subset R$ and $R_k^{\widetilde{u},\widetilde{v}}\subset R$ isomorphic for all $k\in\bbn\cup\{0\}$?
\smallskip

To answer this question, we encounter the following situations.
\smallskip

\noindent\textbf{Case 1:} Let $v\neq\widetilde{v}$. In this case, the subfactors $R_k^{u,v}\subset R$ and $R_k^{\widetilde{u},\widetilde{v}}\subset R$ need not be isomorphic. For example, take $v=F_n$ and $\widetilde{v}$ to be any $n\times n$ complex Hadamard matrix not equal to $F_n$; moreover, choose any $u\in[F_n]$ and $\widetilde{u}$ be any other $n\times n$ complex Hadamard matrix such that the spin model subfactor $R_{\widetilde{u}}\subset R$ is of infinite depth. For example, when $n=4$, choose $\widetilde{u}$ to be a matrix of the following form
\begin{IEEEeqnarray*}{lCl}
\widetilde{u}= \frac{1}{2}\,\left[{\begin{matrix}
1 & 1 & 1 & 1\\
1 & ie^{i\theta} & -1 & -ie^{i\theta}\\
1 & -1 & 1 & -1\\
1 & -ie^{i\theta} & -1 & ie^{i\theta}\\
\end{matrix}}\right]
\end{IEEEeqnarray*}
where $\theta\in(0,\pi)$ is irrational (see \cite{Bur}). It turns out that all the vertex model subfactors $R_k^{u,v}\subset R$ are diagonal subfactors (we shall prove in \Cref{Sec6}, \Cref{diagonal subfactor}), and hence of finite depth. However, since the spin model subfactor $R_{\widetilde{u}}\subset R$ is of infinite depth, and it is an intermediate subfactor to all the vertex model subfactors, it follows that all the vertex model subfactors $R_k^{\widetilde{u},\widetilde{v}}\subset R$ must be of infinite depth by \cite{Bi}. Clearly then, the pairs of subfactors $R_k^{u,v}\subset R$ and $R_k^{\widetilde{u},\widetilde{v}}\subset R$ are nonisomorphic for each $k$.
\smallskip

\noindent\textbf{Case 2:} Let $v=\widetilde{v}$. It turns out that in this case, two subcases arise depending on whether $u\mbox{ and }\widetilde{u}$ are from the same Hadamard class, or they are Hadamard inequivalent (see \Cref{Hadamard equivalence}).
\smallskip

\noindent\textbf{Subcase 2.1:} Let $v=\widetilde{v}$ and $u,\widetilde{u}$ are Hadamard inequivalent. In this case also, $R_k^{u,v}\subset R$ and $R_k^{\widetilde{u},\widetilde{v}}\subset R$ need not be isomorphic. One can consider the example in Case $1$ with identical $u\mbox{ and }\widetilde{u}$, and $v=\widetilde{v}=F_n$. Due to the same reason, one arrives at an exactly similar conclusion as in Case $1$.
\smallskip

\noindent\textbf{Subcase 2.2:} Let $v=\widetilde{v}$ and $u,\widetilde{u}$ be Hadamard equivalent. It turns out that the resulting pairs of vertex model subfactors are always isomorphic in this case. Moreover, we can also characterize when the isomorphism is inner and outer. 

Therefore, the answer to the question at the beginning of the section lies in proving the subcase $2.2$ only. To fix our working module, let $W$ be a complex Hadamard matrix of order $n$ and consider the Hadamard equivalence class $[W]$. By definition of Hadamard equivalence, any element $u\in[W]$ is of the form $u=D_1P_1WP_2D_2$, where $D_1, D_2$ are diagonal unitary matrices and $P_1, P_2$ are permutation matrices (see \Cref{Sec2.1}). Let $u,v\in[W]$ be two distinct complex Hadamard matrices. We obtain the spin model subfactors $R_u\subset R$ and $R_v\subset R$. For two distinct complex Hadamard matrices $u$ and $v$ of order $n$, we define an equivalence relation by the following
\begin{eqnarray}\label{mainequivalence}
u\sim v \mbox{ if }v=uPD,
\end{eqnarray}
where $D$ is a diagonal unitary matrix and $P$ is a permutation matrix each of order $n$. Observe that `$\sim$' is a finer equivalence relation than the Hadamard equivalence described in \Cref{Hadamard equivalence}. Now, recall the following result from \cite{BG1}.
\begin{thm}[Thm. $4.7$, \cite{BG1}]\label{subrelation}
\begin{enumerate}[$(i)$]
\item For distinct $n\times n$ complex Hadamard matrices $u\mbox{ and }v$, the pair of Hadamard subfactors $R_u\subset R\mbox{ and }R_v\subset R$ are distinct $(\mbox{that is, }R_u\neq R_v)$ if and only if $\,u\nsim v$.
\item If two $n\times n$ complex Hadamard matrices $u\mbox{ and }v$ are Hadamard inequivalent, then the corresponding spin model subfactors $R_u\subset R\mbox{ and }R_v\subset R$ are always distinct $(R_u\neq R_v)$.
\end{enumerate}
\end{thm}
Due to \Cref{subrelation}, we see that without loss of generality, we may assume that $u=D_1P_1W$ and $v=D_2P_2W$. This is because the spin model subfactors arising from the complex Hadamard matrices $D_1P_1W\mbox{ and }D_1P_1W\widetilde{P}_1\widetilde{D}_1$ are the same, due to \Cref{mainequivalence} and \Cref{subrelation}. The following result is the complete answer to the subcase $2.2$ above.
 
\begin{thm}\label{conjugate theorem}
Let $u,v$ be Hadamard equivalent but distinct complex Hadamard matrices of the same order, and $R_u, R_v\subset R$ be the corresponding Hadamard/spin model subfactors. Without loss of generality, suppose that $u = D_1 P_1 W$ and $v = D_2 P_2 W$. Let $w$ be any arbitrary complex Hadamard matrix of the same order as $u$ and $v$. Then, we have the following\,: 
\begin{enumerate}[$(i)$]
\item There exists $\theta\in\mathrm{Aut}(R)$ such that the following chain of subfactors
\[
\ldots\subset R_k^{u,w}\subset\ldots\subset R_0^{u,w}\subset R_u\subset R
\]
and
\[
\ldots\subset R_k^{v,w}\subset\ldots\subset R_0^{v,w}\subset R_v\subset R
\]
are isomorphic, that is, $\theta(R_u)=R_v$ and $\theta(R^{u, w}_k)=R^{v, w}_k$ for all $k$.
\item Furthermore, the automorphism $\theta$ is outer if and only if $P_1\neq P_2$, in which case the outer period equals to the order of the permutation $P_2P_1^*$.
\end{enumerate}
\end{thm}

A striking feature in the above theorem is that one single automorphism works for all the vertex model subfactors; in other words, $\theta$ implements a chain isomorphism. The proof involves some effort, and we begin with the following lemma, which has appeared implicitly in the proof of Proposition $1.6$, \cite{Connes}. However, for the sake of completeness, we supply a proof here.

\begin{lmma} \label{morphism}  
Let $R$ be the hyperfinite type-$II_1$ factor and $\{\alpha_n\}_{n \in \mathbb{N}}$ be a sequence of injective homomorphisms from $R$ into $R$ satisfying  
\begin{equation*}
\alpha_k(x)=\alpha_{k+r}(x)
\end{equation*}
for all $x \in A_{2k}$ (see \Cref{not}) and $r \in \bbn$.
Then, the sequence $\{\alpha_n\}_{n \in \mathbb{N}}$ \textit{converges pointwise strongly} to an \textit{injective homomorphism} $\theta$ of $R$.  
\end{lmma}
\begin{prf}
Let $L^2(R, \tau)$ be the GNS Hilbert space with cyclic vector $\Omega$ of the hyperfinite factor $R$ which has unique trace $\tau$. The hyperfinite factor $R$ acts on a dense subspace $R \Omega \subset L^2(R, \tau)$ by left multiplication, and the anti-unitary operator $J$ defined by
\begin{equation}\label{involution}
J x\Omega = x^*\Omega, \quad \text{for all}\, x \in R.
\end{equation} 
satisfies the commutation relation $J R J = R'$, where $R'$ is the commutant of $M$. Note that, from \Cref{involution} we have  
\[
R \Omega = \{JmJ \Omega : m \in R \} \subset L^2(R, \tau).
\]  
\noindent Consider any $x \in R$. Define a sequence $\{x_n\}$ by  
\[
x_n = E^R_{A_{2n}}(x),
\]  
where $E^{R}_{{A}_{2n}}$ denotes the trace-preserving conditional expectation onto $A_{2n}$. Then we know that the sequence $\{x_n\}$ converges to $x$ in the strong operator topology (SOT). Therefore, for any $\epsilon > 0$ and $y \in M$ there exists $n_0 \in \mathbb{N}$ such that
\begin{equation}\label{epsilon bound}
 \|(x_n - x)\Omega \| < \frac{\epsilon}{2 \| JyJ \|},   
\end{equation}
where $\|JyJ\|$ denotes the operator norm of $JyJ: L^2(R, \tau) \to L^2(R, \tau)$. For $n, m > n_0$, it follows from the hypothesis that
\begin{equation}\label{n0}
\alpha_m(x_{n_0}) = \alpha_n(x_{n_0}) = \alpha_{n_0}(x_{n_0})
\end{equation}
Using the commutation property $JRJ = R'$, and \Cref{n0}, for any $n, m \geq n_0$, we obtain the following
\begin{IEEEeqnarray*}{lCl}
    \|(\alpha_{n}(x)-\alpha_{m}(x))JyJ \Omega\|&=&\|JyJ(\alpha_{n}(x)-\alpha_{m}(x))\Omega\|\\
    &\leq&\|JyJ\| \, \|(\alpha_{n}(x)-\alpha_{m}(x)) \Omega\| \nonumber \\
    &=&||JyJ|| \, \|(\alpha_{n}(x)-\alpha_n(x_{n_0})+\alpha_{m}(x_{n_0})-\alpha_{m}(x)) \Omega\|\\
    &\leq&\|JyJ\| \, \|\alpha_{n}(x-x_{n_0})\Omega\|+\|\alpha_{m}(x_{n_0}-x)) \Omega\|\\
    &=&\|JyJ|| \, \|(x-x_{n_0})\Omega\|+\|(x_{n_0}-x)\Omega\| \nonumber \hspace{2mm} (\text{by \Cref{epsilon bound}})\\
    &\leq& \epsilon
\end{IEEEeqnarray*}
hence the sequence $\{\alpha_n(x)\}$ forms a Cauchy sequence in the dense subspace $R \Omega \subset L^2(R, \tau)$ with respect to strong operator topology. Therefore we define the map $\theta: R \to R$ by  
\[
  \theta(x)(JyJ \Omega) = \lim_{n \to \infty} \alpha_n(x)(JyJ \Omega), \quad \forall \, y \in R.  
\]
Since each $\alpha_n$ is an injective homomorphism, we have $\|\alpha_n(x)\| = \|x\|$, ensuring that $\theta(x)$ extends linearly to all of $L^2(R, \tau)$ and by the density of $R \Omega$ in $L^2(R, \tau)$, the map $\theta$ is well-defined and satisfies the pointwise strong convergence property  
\begin{eqnarray}\label{dfn}
\alpha_n(x) \xi \to \theta(x) \xi, \quad \forall \, x \in R \hspace{1mm} \text{and} \hspace{1mm} \xi \in L^2(R, \tau).
\end{eqnarray}
Moreover, since 
\[
\|\alpha_n(x)\| = \|x\|, \quad \|\alpha_n(y)\| = \|y\|,
\]  
it follows easily that $\theta(xy) = \theta(x)\theta(y)$.

Now we prove the injectivity of $\theta$. Since $R$ has a unique trace and each $\alpha_n$ is an injective homomorphisms, we have  
$\|\alpha_n(x)\Omega\| = \|x\Omega\|$ for all $n$. Then, by \Cref{dfn}, we conclude the following\,:
\[
\| \theta(x)\Omega\|=\|x \Omega\|=\sqrt{\tau(x^*x)}.
\]
This immediately implies that $\theta$ is injective.\qed  
\end{prf}

\begin{ppsn}\label{outer automorphism}
Let $u$ be a unitary matrix in $M_n$ and $\{\alpha_k\}$ be the sequence of inner automorphisms of the hyperfinite type $II_1$ factor $R$ defined by
\begin{equation*}
\alpha_{k}(x):=\mathrm{Ad}_{(u \ot u^{(k)})}(x)
\end{equation*}
for all $x \in R$. Then, the sequence $\{\alpha_k\}$ strongly converges to an automorphism $\theta$ of $R$. Furthermore, we have the following:
\begin{enumerate}[$(i)$]
\item If $u\neq I_n$, then $\theta$ is an outer automorphism of $R$.
\item If $m$ is the smallest positive integer such that $u^m=\lambda_0I_n$ for some $\lambda_0\in\bbc$ with modulus one, then $\theta$ is an outer automorphism with outer period  $m$; and if no such $m$ exists, then the outer period of $\theta$ is zero.
\end{enumerate}
\end{ppsn}
\begin{prf}
Let $\tau$ is a unique trace on hyperfinite $II_1$ factor $R=\overline{\cup_j A_{2j}}^{\,\mbox{sot}}$. Define $\|x\|_2=\sqrt{\tau(x^*x)}$ for all $x \in R$. Observe that for any fixed $k\in\bbn$ and $x\in A_{2k}$, we have 
\begin{equation}\label{constant sequence}
\alpha_{k+r}(x)=\alpha_k(x)
\end{equation}
for all $r\in\bbn$, and therefore from \Cref{morphism}, we conclude that there exists an injective homomorphism $\theta$ such that the sequence $\{\alpha_n(x)\}$ converges pointwisely strongly to $\theta(x)$ for all $x \in R$. By hypothesis and \Cref{constant sequence}, we get $\alpha_{k}(A_{2k})=\alpha_{k+r}(A_{2k})=A_{2k}$ for all $r\in\bbn$, and hence, it follows that $\theta(A_{2k})=A_{2k}$ for all $k\in\bbn$. Therefore $\cup_kA_{2k}\subset\theta(R)$, and by the continuity of the injective homomorphism $\theta$, we conclude that $\theta$ is an onto map. Hence, $\theta$ is an automorphism of $R$.

Now, let $u\neq I_n$. We show that $\theta$ is an outer automorphism of $R$. Consider the shift endomorphism $\pi$ on $R$ defined by 
\[
\pi(x_1 \ot x_2 \ot \cdots \ot x_n)=(x_1 \ot x_2 \ot \cdots \ot x_n \ot I_n).
\]
For any $y \in M_n$, observe that $\{\pi^{k}(y)\}_{k \in \bbn}$ is a central sequence in $R$. That is, for any $x\in R$, we have the following
\[
||x \pi^{k}(y)-\pi^{k}(y)x||_{2} \rightarrow 0\hspace{4mm} \text{when}  \hspace{1mm } k\rightarrow \infty.
\]
Therefore, if $\beta$ is an inner automorphism of $R$, then 
\[
||\beta( \pi^{k}(y))-\pi^{k}(y)||_{2} \rightarrow 0, \hspace{4mm} \text{when}  \hspace{1mm } k\rightarrow \infty.
\]
Take any unitary matrix $w\in M_n$ such that $\text{Ad}_{w}(u)\notin\dt$, the diagonal subalgebra of $M_n$. Then for $k\geq 2$, we have the following
\begin{eqnarray}\label{second normed equation}
||\theta(\pi^{k}(w^{*}\mathscr{D}_1 w))-\pi^{k}(w^{*}\mathscr{D}_1 w)||_2 &=& ||\pi^{k}(uw^{*} \mathscr{D}_1 w u^{*})- \pi^{k}(w^{*} \mathscr{D}_1 w)||_2 \nonumber \\
&=& ||uw^{*} \mathscr{D}_1 w u^{*}- w^{*} \mathscr{D}_1 w||_2.
\end{eqnarray}
Observe that $||uw^{*} \mathscr{D}_1 w u^{*}- w^{*} \mathscr{D}_1 w||_2\neq 0$. This follows because on contrary, if $||uw^{*} \mathscr{D}_1 w u^{*}- w^{*} \mathscr{D}_1 w||_2=0$, then $\text{Ad}_{wuw^{*}}{(\mathscr{D}_1)}={\mathscr{D}_1}$, which further implies that $\text{Ad}_{wuw^{*}}(D)=D$ for all $D \in \dt$, which contradicts the assumption that $\text{Ad}_{w}(u)\notin\dt$. Therefore, from \Cref{second normed equation}, we obtain the following
\[
||\theta(\pi^{k}(w^{*}\mathscr{D}_1 w))-\pi^{k}(w^{*}\mathscr{D}_1 w)||_2 \not \rightarrow 0 \hspace{4mm} \text{whenever}  \hspace{1mm } k\rightarrow \infty.
\]
This concludes that $\theta$ is an outer automorphism of $R$, which finishes the proof of part $(i)$. 

For part $(ii)$, first observe that for $x \in A_{2k}$, we have $\alpha_k(x)\in A_{2k}$, and hence from \Cref{constant sequence} it follows that 
\begin{equation}\label{power}
\theta^{m}(x)=\lim_k\alpha^{m}_{k}(x)
\end{equation}
for all $m\in\bbn$. Since $\{\alpha^{m}_k\}$ is the sequence of inner automorphisms satisfying $\alpha^{m}_{k+r}(x)=\alpha^{m}_{k}(x)$ for all $x \in A_{2k}$ and $r\in\bbn$, by using \Cref{power} and \Cref{morphism}, it finally follows that
\begin{equation}\label{gsk}
\theta^{m}(x)=\text{lim}_{k \rightarrow  \infty}\,\alpha^{m}_{k}(x) \hspace{3mm} \text{for all} \hspace{1mm} x \in R=\overline{\cup_kA_{2k}}.
\end{equation}
Let $m\in\{1,2,...,n-1\}$ and define $v=u^{m}$. By the given hypothesis, $v\neq I_n$. Observe that $\alpha^{m}_{k}=\text{Ad}_{v^{(k+1)}}$, and therefore from \Cref{gsk} and part $(i)$ of the statement, we conclude that $\theta^{m}$ is an outer automorphism of $R$. By our hypothesis, it follows that $\alpha^{n}_k=\mbox{id}_{R}$, where $\mbox{id}_{R}$ denotes the identity automorphism of $R$. Therefore, from \Cref{gsk} it follows that $\theta^n=\mbox{id}_R$, and this completes the proof of part $(ii)$.\qed
\end{prf}
\smallskip

\noindent\textbf{Proof of \Cref{conjugate theorem}:} Take $P=P_2 P^{*}_1$ and then putting $u=P$ in \Cref{outer automorphism}, we obtain $\theta\in\mathrm{Aut}(R)$ satisfying
\begin{equation}\label{finite stage}
\theta(x)=\alpha_{k}(x)\quad\mbox{for all }\,x \in A_{2k}
\end{equation}
where $\alpha_k=\mathrm{Ad}_{P^{(k+1)}}$. Since $P$ is a permutation, \Cref{finite stage} gives us the following\,:
\begin{equation}\label{gsk1}
\theta(e_{2i+2})=e_{2i+2}\,,\,\,\theta(e_{2i+1})=e_{2i+1}\,,\,\,\theta(u D u^{*})=\text{Ad}_{P}(u D u^{*}),
\end{equation}
where $e_{2i+2}=e_2 \ot I^{(i)}_n$ and $e_{2i+1}=e_1 \ot I^{(i)}_n \in A_{2i}$ are the Jone's projections and $D\in\Delta_n$. Now, define $\theta_1=\text{Ad}_{\{D_2PD^{*}_1 P^{*}\}} \circ \theta$ and using \Cref{gsk1} observe that
\begin{equation}\label{gsk2}
\theta_1(e_{2i+2})=e_{2i+2}\,,\,\,\theta_1(e_{2i+1})=e_{2i+1}\,,\,\,\theta_1(u D u^{*})=v D v^{*}
\end{equation}
for all $i \in \bbn \cup \{0\}$. We claim that
\begin{IEEEeqnarray}{lCl}\label{guru0}
\theta_1(\text{Ad}_{u_{k}}(x))=\text{Ad}_{v_{k}}(x)
\end{IEEEeqnarray}
for all $x \in A_{k-1}$ and $k \in \bbn$. First observe that from \Cref{not}, we have $e_{2k+1}=e_1 \ot I^{(k)}_n$, $B^{u}_{2k}=\text{Ad}_{u_{2k}}(\dt \ot M^{(k)}_n)$, and similarly, $B^{v}_{2k}=\text{Ad}_{v_{2k}}(\dt \ot M^{(k)}_n)$. Assume that the claim is true for some even step $k=2m$, that is, $\theta_1(\text{Ad}_{u_{2m}}(x))=\text{Ad}_{v_{2m}}(x)$ for all $x \in A_{2m-1}=\dt \ot M^{(m)}_n$. From \Cref{tower unitary}, it follows that
\[
\text{Ad}_{u_{2m+1}}(I_n\otimes x)=\text{Ad}_{u_{2m}}(x)\quad\mbox{ and }\quad\text{Ad}_{v_{2m+1}}(I_n\otimes y)=\text{Ad}_{v_{2m}}(y)
\]
for all $x,y\in\dt \ot M^{(m)}_n$. Therefore, we get the following\,:
\begin{IEEEeqnarray*}{lCl}
& & \theta_{1}(\text{Ad}_{u_{2m+1}}((I_n\ot x)e_{2m+1}(I_n\ot y)))\\
&=&\theta_1\big(\text{Ad}_{u_{2m+1}}(I_n\ot x) \text{Ad}_{u_{2m+1}}(e_{2m+1}) \text{Ad}_{u_{2m+1}}(I_n\ot y)\big)\\
&=&\text{Ad}_{(I_n \ot v_{2m})}(x)) \theta_1(\text{Ad}_{u_{2m+1}}(e_{2m+1})) (\text{Ad}_{(I_n \ot v_{2m})}(y)\qquad\hfill{(\text{apply \Cref{unitary projection relation}})}\\
&=&\text{Ad}_{(I_n \ot v_{2m})}(x)) \theta_1(e_{2m+2}) (\text{Ad}_{(I_n \ot v_{2m})}(y)\hfill{(\text{apply \Cref{gsk2}})}\\
&=&\text{Ad}_{(I_n \ot v_{2m})}(x)) (e_{2m+2}) (\text{Ad}_{(I_n \ot v_{2m})}(y)\hfill{(\text{apply \Cref{unitary projection relation}})}\\
&=&\text{Ad}_{(v_{2m+1})}(x)) \text{Ad}_{(v_{2m+1})}(e_{2m+1}) (\text{Ad}_{(v_{2m+1})}( y)\hfill{(\text{apply \Cref{unitary projection relation}})}\\
&=&\text{Ad}_{v_{2m+1}}(x \,e_{2m+1}\, y).
\end{IEEEeqnarray*}
Hence, we conclude that $\theta_1(\text{Ad}_{u_{2m+1}}(x))=\text{Ad}_{v_{2m+1}}(x)$ for every $x \in M_n \ot M^{(m)}_n$. Following a similar argument, the even step of the claim can also be proven by assuming its previous odd step. This concludes the proof of the claim in \Cref{guru0}.

Now, we show the following  
\[
\theta_1(R_u) = R_v \quad \text{and} \quad \theta_1(R_{u,w}^k) = R_{v,w}^k 
\]  
for all $k \in \mathbb{N} \cup \{0\}$. From the construction of the spin model subfactors $R_u \subset R$ and $R_v \subset R$, we have  
\[
R_u=\overline{\bigvee\{u \dt u^{*}, e_{i+1} : i \in \bbn\}}^{\text{sot}}
\]
and
\[
R_v=\overline{\bigvee\{v \dt v^{*}, e_{i+1} : i \in \bbn\}}^{\text{sot}}.
\]
By the continuity of automorphism $\theta_1 \in \text{Aut}(R)$, and by using  \Cref{gsk2}, we conclude that $\theta_1(R_u) = R_v$. Recall from \Cref{middle vertex subfactor} the following\,:  
\begin{eqnarray*}
R_{k}^{u,w}=\overline{\bigvee \{N^{u,w}_{2k}, \,e_{2k+2+i} \,: \, i\in \bbn\}}^{\text{sot}} \nonumber\\
R_{k}^{v,w}=\overline{\bigvee \{N^{v,w}_{2k}, \,e_{2k+2+i} \,: \, i\in \bbn\}}^{\text{sot}}
\end{eqnarray*}
where
\begin{eqnarray}
N^{u, w}_{2k} &=& \text{Ad}_{u_{2k+1}(I_n \otimes w_{2k})}(I_n \otimes \Delta_n \otimes I_n^{(k)}) \nonumber\\
N^{v, w}_{2k} &=& \text{Ad}_{v_{2k+1}(I_n \otimes w_{2k})}(I_n \otimes \Delta_n \otimes I_n^{(k)})  \, \,(\text{see \Cref{notation two}}).\nonumber
\end{eqnarray}
Since $\text{Ad}_{(I_n \otimes w_{2k})}(I_n \otimes \Delta_n \otimes I_n^{(k)}) \subset A_{2k}$, using \Cref{gsk2} it follows that  
\begin{equation*}
\theta_1(N^{u,w}_{2k}) = N^{v,w}_{2k},
\end{equation*}
and then using the continuity of the automorphism $\theta_1 \in \text{Aut}(R)$, we finally conclude that  
\[
\theta_1(R^{u,w}_k) = R^{v,w}_k.
\]
This completes the proof of the first assertion.
\smallskip

We divide the proof of the second assertion into two cases\,:

\noindent\textbf{Case 1.} $P_1=P_2$.

Since $P_2P^{*}_1=I_n$, by definition of $\alpha_k$ and $\theta$, it follows that $\theta=id_{R}$. Hence, $\theta_1=\text{Ad}_{D_2PD^{*}_1 P^{*}}$ (as $\theta_1=\mathrm{Ad}_{D_2PD^{*}_1 P^{*}} \circ \theta$). This concludes that $R_{v}=\mathrm{Ad}_{D_2PD^{*}_1 P^{*}}(R_u)$ and $R^{v,w}_k=\text{Ad}_{D_2PD^{*}_1 P^{*}}(R^{u,w}_k)$ for all $k \in \bbn \cup \{0\}$.
\smallskip

\noindent\textbf{Case 2.} $P_1\neq P_2$.

Since $P_2P^{*}_1\not = I_n$, from \Cref{outer automorphism} it follows that $\theta$ is an outer automorphism of $R$ with outer period $r$, where $r$ is the least positive integer such that $(P_2P^{*}_1)^{r}=I_n$. Now, for any $m \in \bbn$ putting $\xi=D_2PD^{*}_1 P^{*}$, a straightforward computation gives us $\theta^{m}_1=\text{Ad}_\xi\mathrm{Ad}_{\theta(\xi)}\ldots\mathrm{Ad}_{\theta^{m-1}(\xi)}\circ \theta^{m}$. Therefore, $\theta_1$ is also an outer automorphism with outer period $r$, where $\theta^{r}_1=\mathrm{Ad}_{\xi\theta(\xi)\ldots\theta^{r-1}(\xi)}$, which completes the proof.\qed
\smallskip

\noindent\textbf{Conclusion\,:} Among the pairs $(u,v),\,(\widetilde{u},\widetilde{v})$ as input data, if we are interested in isomorphic pairs of vertex model subfactors $(R_k^{u,v}\subset R) \cong (R_k^{\widetilde{u},\widetilde{v}}\subset R)$ for each $k\in\bbn\cup\{0\}$, as illustrated in \Cref{fig2}, it turns out that $v=\widetilde{v}$ is an `optimal' choice in some sense. To be precise, if $v\neq\widetilde{v}$, we have a counterexample that the resulting pairs of vertex model subfactors need not be isomorphic. This is described in Case $1$ at the beginning of the section. Secondly, even if $v=\widetilde{v}$, it further turns out that if $u\mbox{ and }\widetilde{u}$ are Hadamard inequivalent, then the resulting pairs of vertex model subfactors need not be isomorphic. This is described in Subcase $2.1$. Finally, we are left with only one scenario (Subcase 2.2): when $v=\widetilde{v}$ and $u,\widetilde{u}$ are Hadamard equivalent. In this case, \Cref{conjugate theorem} provides a complete answer.

We propose the following important questions.
\smallskip

\noindent\textbf{Open question 1\,:} Let $v\neq\widetilde{v}$, or $v=\widetilde{v}$ but $u\mbox{ and }\widetilde{u}$ are Hadamard inequivalent. Are the pairs of vertex model subfactors $R^{u,v}_k\subset  R$ and $R^{\widetilde{u},\widetilde{v}}_k\subset R$ in \Cref{fig2} always non-isomorphic? In other words, if the pairs of vertex model subfactors $R^{u,v}_k\subset  R$ and $R^{\widetilde{u},\widetilde{v}}_k\subset R$ are isomorphic, is it necessary that $v=\widetilde{v}$ and $u\mbox{ and }\widetilde{u}$ are Hadamard equivalent?
\smallskip

Observe that the `sufficiency' part is established in \Cref{conjugate theorem}.
\smallskip

\noindent\textbf{Open question 2\,:} Consider \Cref{fig2}. Can we compute the Pimsner-Popa probabilistic constant $\lambda(R_k^{u,v},R_k^{\widetilde{u},\widetilde{v}})$ \cite{PP} and the Connes-St{\o}rmer relative entropy $H(R_k^{u,v}|R_k^{\widetilde{u},\widetilde{v}})$ \cite{CS} for the pairs of subfactors $R_k^{u,v}\subset R\mbox{ and }R_k^{\widetilde{u},\widetilde{v}}\subset R$ for each $k\in\bbn\cup\{0\}$?
\smallskip

\noindent\textbf{Open question 3\,:} What can we say about $R_k^{u,v}\cap R_k^{\widetilde{u},\widetilde{v}}$? Are these (finite-index) subfactors of $R$?


\newsection{The special case of DFT matrices}\label{Sec6}

In this section, we deal with the special case where our input data comes from the Hadamard equivalence class of the Fourier matrix $F_n$, a special type of complex Hadamard matrix. This situation is of utmost importance and interesting, as it helps us to construct many counterexamples. We completely characterize the vertex model subfactors constructed in \Cref{theorem2} for this special situation.

Start with the input data $(u,F_n)$, where $F_n$ is the Fourier matrix of order $n$ and $u\in[F_n]$, the Hadamard equivalence class of $F_n$. We first show that all the vertex model subfactors $R_k^{u,F_n}\subset R$ constructed in \Cref{theorem2} are diagonal, for each $k\in\bbn\cup\{0\}$. Then, we show that this choice of input data is in some sense `optimal' if we restrict outcome as diagonal subfactors. More precisely, if we start with the input data $(u,v)$ of complex Hadamard matrices of the same order, in any other situations except $v=F_n$ and $u\in[F_n]$, the resulting vertex model subfactors need not be diagonal; whereas, in the case of $v=F_n$ and $u\in[F_n]$ they are always so.

Recall from \Cref{theorem2}, the following nested sequence of vertex model subfactors\,:
\[
\cdots\subset R_{k+1}^{u, F_n} \subset R_{k}^{u, F_n} \subset \cdots \subset R_{1}^{u, F_n} \subset R_{0}^{u, F_n} \subset R
\]
arising from the pair $(u,F_n)$, where $u\in[F_n]$. It turns out that one has to focus only on the case $(F_n,F_n)$, thanks to \Cref{conjugate theorem}. So, it is good enough to characterize the following nested sequence of vertex model subfactors
\[
\cdots\subset R_{k+1}^{F_n, F_n} \subset R_{k}^{F_n, F_n} \subset \cdots \subset R_{1}^{F_n, F_n} \subset R_{0}^{F_n, F_n} \subset R
\]
where  
\begin{eqnarray}\label{vertex plus fourier}
R_{k}^{F_n,F_n}=\overline{\bigvee \{N^{F_n, F_n}_{k}, \,e_{2(k+1)+i} \,: \, i\in \bbn\}}^{sot} 
\end{eqnarray}
(see \Cref{middle vertex subfactor}). Our main result is the following.
\begin{thm}\label{diagonal subfactor}
Let $F_n$ be the Fourier matrix of order $n$. There is an outer automorphism $\theta\in\mathrm{Out}(R)$ with outer period $n$ such that the vertex model subfactor $R_k^{F_n,F_n}\subset R$ constructed in \Cref{theorem2} is isomorphic to the following diagonal subfactor
\[
\mathscr{R}_k:=\bigg{\{} \displaystyle{\sum_{j=1}^{n^{k+1}} E_{jj} \ot \theta^{\,r_j}(x) : \, x \in R,\,r_j\in\{0,\ldots,n-1\}},\,\theta^n=\mathrm{id}\bigg{\}} \subset M^{(k+1)}_n\ot R
\]
for each $k\geq 0$. Furthermore, their relative commutants are given by the following\,:
\[
\big(R^{F_n, F_n}_k\big)^\prime\cap R=\bbc^{n^{k+1}}.
\]
\end{thm}

\begin{rmrk}\rm
\begin{enumerate}[$(i)$]
\item A striking feature of the above theorem is that we have constructed an outer automorphism $\theta$ that works for all the vertex model subfactors $R_k^{F_n,F_n}\subset R,\,k\geq 0$. That is, $\theta$ is independent of $k$, and its outer period is exactly equal to the order of the input matrix.
\item The above theorem gives the proof of the example mentioned in Case $1$, \Cref{Sec5}.
\end{enumerate}
\end{rmrk}

We start by proving a couple of lemmas. First, we fix the following notations throughout this section.

\begin{notation}\label{notation3}\rm
Let $\omega$ denote the primitive $n$-th root of unity. Fix the following diagonal unitary and permutation matrices in $M_n:$
\begin{align*}
       \mathscr{D}_{0}:=I_n \,,\quad \quad \quad \quad \quad  \quad& \sigma_{1}:=E_{12}+E_{23}+ \cdots E_{(n-1)n}+E_{n1}\,,\\ 
        \mathscr{D}_{1}:=\mbox{diag} \{1,\omega, \omega^{2}, \cdots ,\omega^{(n-1)}\} \,,\quad \quad \quad & \sigma_{k}:=\sigma^{k}_{1}\,.\\
        \hspace{4mm}\mathscr{D}_{k}:=\mathscr{D}^{k}_{1}=\mbox{diag} \{1,\omega^{k}, \omega^{2k}, \cdots ,\omega^{(n-1)k}\} \,.\quad&      
\end{align*}
Observe that $\sigma_1$ and $\mathscr{D}_1$ denote the adjoint of the shift matrix and the clock matrix respectively.
\end{notation}

\begin{lmma}\label{diagonal permutation relation}
For $k=0,1,\cdots,(n-1)$, we have the following identities
\begin{eqnarray*}
\mathscr{D}_{1}F_{n}=F_{n}\sigma_{n-1}\quad,\quad\quad\sigma_{1}F_n=F_n\mathscr{D}_{1}\,,\nonumber\\
\mathscr{D}_{k}F_{n}=F_{n}\sigma^{k}_{n-1}\quad,\quad\quad\sigma_{k}F_n=F_n\mathscr{D}_{k}\,.
\end{eqnarray*}
\end{lmma}
\begin{prf}
These are simple verification, and we omit the details.\qed
\end{prf}

For the case of Fourier matrix the $F_n$, we have the following unitary matrices\,:
\begin{eqnarray}\label{Fourier unitaries}
(F_n)_{2k+1}=(I_n \otimes (F_n)_{2k})(D_{F_n} \otimes I_n^{(k)})  \quad,\quad (F_n)_{2k}=(F_n)_{2k-1}(F_n \otimes I_n^{(k)}),
\end{eqnarray}
from \Cref{tower unitary}, where $D_{F_n}:=\text{bl-diag}\{I_n, \,\mathscr{D}_{n-1}, \cdots , \,\mathscr{D}_2, \,\mathscr{D}_1\}$ is a block-diagonal matrix.
    
\begin{lmma}\label{first lemma}
For $k \in \bbn \cup \{0\}$, we have the following identity
\begin{equation*}
 \mathrm{Ad}_{(I_n \ot (F_n)_{2k})}(I_n \ot \mathscr{D}_1 \ot I^{(k)}_{n})=I_n \ot \sigma_1^{(k+1)}.
\end{equation*}
where $(F_n)_{2k}$'s are the unitary matrices as defined in \Cref{Fourier unitaries}. 
\end{lmma}
\begin{prf}
From \Cref{diagonal permutation relation}, we have $F_n \mathscr{D}_j=\sigma_{n-j} F_n$, and therefore it follows that 
\begin{eqnarray}\label{block permutation matrix}
 (F_n)_1&=&W_2 (I_n \ot F_n),
\end{eqnarray}
where $(F_n)_1=(I_n \ot F_n)D_{F_n}$ (see \Cref{Fourier unitaries}) and $W_2=\text{bl-diag}\{I_n, \sigma_{n-1}, \cdots, \sigma_2, \sigma_1\}$. We prove by induction on $k$. From \Cref{diagonal permutation relation}, it easily follows that 
\begin{equation*}
\text{Ad}_{(I_n \ot F_n)}(I_n \ot \mathscr{D}_1)=I_n \ot \sigma_1,
\end{equation*}
which concludes the base step $k=0$. Now, assume that the statement is true for the $k$-th step. From \Cref{Fourier unitaries}, one easily get that $(F_n)_{2k+2}=(I^{(k)}_n \ot (F_n)_1)((F_n)_{2k} \ot I_n)$. We use this to prove the statement for the $(k+1)$-th step.
\begin{IEEEeqnarray*}{lCl}
     & &\text{Ad}_{(I_n \ot (F_n)_{2k+2})}(I_n \ot \mathscr{D}_1 \ot I^{(k+1)}_{n})\\
     &=&\text{Ad}_{(I^{(k+1)}_n \ot (F_n)_1)( I_n \ot (F_n)_{2k} \ot I_n)}(I_n \ot \mathscr{D}_1 \ot I^{(k+1)}_n)\\ 
     &=&\text{Ad}_{(I^{(k+1)}_n \ot (F_n)_1)}\text{Ad}_{(I_n \ot (F_n)_{2k} \ot I_n)}(I_n \ot \mathscr{D}_1 \ot I^{(k+1)}_n)\\
     &=&\text{Ad}_{(I^{(k+1)}_n \ot (F_n)_1)}\big{(}\text{Ad}_{\{I_n \ot (F_n)_{2k}\}}(I_n \ot \mathscr{D}_1 \ot I^{(k)}_n)\big{)}\ot I_n\qquad\hfill{(\text{by induction hypothesis})}\\
     &=&\text{Ad}_{(I^{(k+1)}_n \ot (F_n)_1)} ((I_n \ot \sigma^{(k+1)}_1)\ot I_n)\\
     &=&I_n \ot \sigma^{(k)}_1 \ot \text{Ad}_{(F_n)_1} (\sigma_1 \ot I_n) \qquad\hfill{(\Cref{block permutation matrix})}\\
     &=&I_n \ot \sigma^{(k)}_1 \ot \text{Ad}_{W_2} (\sigma_1 \ot I_n)\\
     &=&I_n \ot \sigma^{(k)}_1  \ot \sigma_1 \ot \sigma_1.
\end{IEEEeqnarray*}
This completes the induction.\qed
\end{prf}

From \Cref{notation two}, we have the following\,:
\begin{eqnarray}\label{fourier algebra}
N^{F_n, F_n}_{2k}=\text{Ad}_{(F_n)_{2k+1}(I_n \otimes (F_n)_{2k})}(I_n \otimes \dt \otimes I_n^{(k)}).
\end{eqnarray}

\begin{lmma}\label{diagonal von Neumann algebra}
We have the following identity
\[
\mathrm{Ad}_{(F_n)_{2k+1}(I_n \ot (F_n)_{2k})}(I_n \ot \mathscr{D}_1 \ot I^{(k)}_{n})=\displaystyle{\sum_{i, i_1,\cdots,i_k=1}^{n}\mathrm{Ad}_{\sigma_{r_i}} (\mathscr{D}_1) \ot E_{i_1i_1}\ot \cdots \ot E_{i_{k}i_{k}} \ot E_{ii}}
\]
where $\mathscr{D}_1$ is defined in \Cref{notation3}, $r_{i} \in \{0,1,2,\dots,n-1\}$, and $(F_n)_{2k}, (F_n)_{2k+1}$ are the unitary matrices as defined in \Cref{Fourier unitaries}. 
\end{lmma}
\begin{prf}
From \Cref{first lemma}, we have $\text{Ad}_{(I_n \ot (F_n)_{2k})}(I_n \ot \mathscr{D}_1 \ot I^{(k)}_{n})=I_n \ot \sigma_1^{(k+1)}$, and thus we obtain the following\,:
\begin{eqnarray*}
\text{Ad}_{(F_n)_{2k+1}(I_n \ot (F_n)_{2k})}(I_n \ot \mathscr{D}_1 \ot I^{(k)}_{n})&=&\text{Ad}_{(F_n)_{2k+1}}(I_n \ot \sigma^{(k+1)}_1).
\end{eqnarray*}
Therefore, it suffices to establish the following identity\,:
\begin{equation}\label{new equation}
\text{Ad}_{(F_n)_{2k+1}}(I_n \ot \sigma^{(k+1)}_1)=\displaystyle{\sum_{i,i_1,\ldots,i_k=1}^{n}\text{Ad}_{\sigma_{r_i}} (\mathscr{D}_1) \ot E_{i_1i_1}  \ot \cdots \ot E_{i_{k}i_{k}} \ot E_{ii}}\,,
\end{equation}
where $r_i\in\{0,\ldots,n-1\}$. We prove \Cref{new equation} by induction. Since $F_n \sigma_1 F^{*}_n=\mathscr{D}_{n-1}$, we obtain the following
\begin{IEEEeqnarray}{lCl}\label{basis step}
\text{Ad}_{(F_n)_1}(I_n \ot \sigma_1)&=&\text{Ad}_{W_2}(I_n \ot \mathscr{D}_{n-1}) \qquad\hfill{(\text{as} \,(F_n)_1=W_2(I_n \ot F_n))}\nonumber\\ 
&=& \sum_{k=1}^{n} E_{ii} \ot \sigma^{*}_{i-1}\mathscr{D}_{n-1} \sigma_{i-1}\nonumber\\
&=& \sum_{i=1}^{n} \text{Ad}_{\sigma_{i-1}}( \mathscr{D}_1) \ot E_{ii}.
\end{IEEEeqnarray}
Therefore, \Cref{new equation} holds for the base step (i.e., for $k=0$). Assume that it holds for the $k^{th}$ step, that is,
\begin{equation}\label{hypothesis}
\text{Ad}_{(F_n)_{2k+1}}(I_n \ot \sigma^{(k+1)}_1)=\displaystyle{\sum_{i,i_1,\cdots,i_k=1}^{n}\text{Ad}_{\sigma_{r_i}} (\mathscr{D}_1) \ot E_{i_1i_1}  \ot \cdots \ot E_{i_{k}i_{k}} \ot E_{ii}}
\end{equation}
where $r_i\in\{0,1,\dots,n-1\}$. To prove the $k+1$-th step, first note that from \Cref{Fourier unitaries} we get $(F_n)_{2k+3}=(I_n^{(k+1)} \ot (F_n)_1)((F_n)_{2k+1} \ot I_n)$, and therefore,
\begin{IEEEeqnarray}{lCl}\label{final equation}
    &&\text{Ad}_{(F_n)_{2k+3}}(I_n \ot \sigma^{(k+2)}_1)\nonumber\\
     &=&\text{Ad}_{(I_n^{(k+1)} \ot (F_n)_1)}\text{Ad}_{((F_n)_{2k+1} \ot I_n)}(I_n \ot \sigma^{(k+2)}_1)\nonumber\\
    &=&\text{Ad}_{(I_n^{(k+1)} \ot (F_n)_1)}\bigg{(}(\text{Ad}_{(F_n)_{2k+1}\}}(I_n \ot \sigma^{(k)}_1))\ot \sigma_1\bigg{)}\quad\hfill{(\text{by \,\Cref{hypothesis}})}\nonumber\\
     &=&\text{Ad}_{(I_n^{(k+1)} \ot (F_n)_1)}\bigg{(}\displaystyle{ \sum_{j,i_1,\cdots,i_k=1}^{n}\text{Ad}_{\sigma_{r_j}} (\mathscr{D}_1) \ot E_{i_1i_1}  \ot \cdots \ot E_{i_{k}i_{k}} \ot E_{jj} \ot \sigma_1}\bigg{)}\nonumber\\
     &=&\bigg{(}\displaystyle{\sum_{j} \sum_{i_1,\cdots,i_k=1}^{n}\text{Ad}_{\sigma_{r_j}} (\mathscr{D}_1) \ot E_{i_1i_1}  \ot \cdots \ot E_{i_{k}i_{k}} \ot \text{Ad}_{ (F_n)_1} (\sum_{i=1}^{n} E_{jj} \ot \sigma_1)}\bigg{)}\hfill{(\text{by \Cref{basis step}})}\nonumber\\
    &=&\displaystyle{\sum_{j}\sum_{i, \,i_1,\cdots,i_k=1}^{n}\text{Ad}_{\sigma_{r_j}} (\mathscr{D}_1) \ot E_{i_1i_1}  \ot \cdots \ot E_{i_{k}i_{k}} \ot {\text{Ad}_{\sigma_{i-1}}(\mathscr{D}_{1})}\ot E_{ii}}
\end{IEEEeqnarray}
where $r_j\in\{0,1,\dots,n-1\}$. It is easy to verify that 
\begin{IEEEeqnarray*}{lCl}
    \text{Ad}_{\sigma_{i-1}}(\mathscr{D}_1)&=&\text{diag}\{\omega^{i},\omega^{i+1},\cdots,\omega^{(n-1)}, \,1,\, \omega, \cdots,w^{(i-1)}\}\\
    &=&\omega ^{i} \mathscr{D}_1.
\end{IEEEeqnarray*}
Hence, applying this fact to \Cref{final equation} deduces to
\[
\displaystyle{\sum_{i, i_1,\cdots,i_k,i_{k+1}=1}^{n}\text{Ad}_{\sigma_{r_i}} (\mathscr{D}_1) \ot E_{i_1i_1}  \ot \cdots \ot E_{i_{k}i_{k}} \ot E_{i_{(k+1)}i_{(k+1)}}\ot E_{ii}}
\]
where $r_i \in\{0,1,2,\cdots,n-1\}$, and the proof is concluded.\qed
\end{prf}

\begin{lmma}\label{make it diagonal}
Let $k\in\bbn\cup\{0\}$. Consider the unitary elements $W_k:=\displaystyle{\sum_{i,j=1}^{n^{k+1}} E_{ij} \otimes E_{ji}}$ and $V_k:=I_n^{(k+1)}\otimes\big(\sum_{i,j=1}^n E_{ij} \otimes I_n^{(k)} \otimes E_{ji}\big)$ in the hyperfinite factor $M_n^{(k+1)} \otimes R$. Define the injective homomorphisms $\beta_k: R \to M_n^{(k+1)} \otimes R$ and $\Gamma_k: M_n^{(k+1)} \otimes R \to M_n^{(k+1)} \otimes R$ given by
\begin{equation}
\Gamma_k(x)=\mathrm{Ad}_{V_kW_k}(x), \quad \beta_k(y) = I_n^{(k+1)} \otimes y
\end{equation}
where $x \in M_n^{(k+1)} \otimes R$ and $y \in R$. Then, the composition $\Gamma_k \circ \beta_k: R \to M_n^{(k+1)} \otimes R$ is an isomorphism.    
\end{lmma}
\begin{prf}
Fix $k\in\bbn\cup\{0\}$. We only need to prove the surjectivity of $\Gamma_k \circ \beta_k$. We claim the following\,:
\begin{eqnarray}\label{guru1}
\Gamma_k \circ \beta_k(M_n^{(k+s+1)})=M^{(k+1)}_n \ot M^{(s)}_n
\end{eqnarray}
for all $s\in\bbn$. Take any element $y_s = x_{k+s+1} \otimes \cdots \otimes x_k \otimes \cdots \otimes x_2 \otimes x_1 \in M_n^{(k+s+1)} \subseteq R$, where $s\in\bbn$. By hypothesis, we have $W_k=\displaystyle{\sum_{i,j}^{n^{k+1}}E_{ij} \ot E_{ji}} \in M^{(k+1)}_n \ot R$, which we can rewrite as $W_k=\sum_{i,j=1}^{n^{k+1}} E_{ij}\otimes I_n^{(k+s)}\otimes E_{ji}$. Hence, 
\begin{eqnarray}\label{apply second conjugation}
\text{Ad}_{W_k}(I^{(k+1)}_n\ot y_s) &=& x_{k+1}\ot x_{k}\ot\cdots\ot x_2\ot x_1\ot x_{k+s+1}\ot\cdots\ot x_{k+2}\ot I^{(k+1)}_n. 
\end{eqnarray}
Rewriting $V_k=\sum_{i,j}^n I_n^{(k+1)} \otimes I_n^{(s)} \otimes E_{ij} \otimes I_n^{(k)} \otimes E_{ji}$ for all $s\in \bbn$, and using \Cref{apply second conjugation}, we obtain the following 
\begin{IEEEeqnarray}{lCl}\label{guru20}
\Gamma_{k}\circ \beta_{k}(y_s) &=&\text{Ad}_{V_k} \text{Ad}_{W_k}(I^{(k+1)}_n \ot y_s)\nonumber\\
&=&\text{Ad}_{V_k}(x_{k+1} \ot x_{k} \ot \cdots \ot x_2 \ot x_1 \ot x_{k+s+1} \ot \cdots \ot x_{k+2} \ot I^{(k+1)}_n)\nonumber\\
&=&x_{k+1} \ot x_{k} \ot \cdots \ot x_2 \ot x_1 \ot x_{k+s+1} \ot \cdots \ot I_n \ot I^{(k)}_n \ot x_{k+2}\nonumber\\
&=&x_{k+1} \ot x_{k} \ot \cdots \ot x_2 \ot x_1 \ot x_{k+s+1}\ot \cdots  \ot x_{k+2},
\end{IEEEeqnarray}
which concludes the proof of the claim \Cref{guru1}. Therefore, it follows that $\Gamma_k \circ \beta_k (\cup_{j} A_{2j})=\cup_{i} (M^{(k+1)}_n \ot A_{2i})$ (see \Cref{not}). Since the injective homomorphism $\Gamma_k \circ \beta_k$ is continuous with respect to the strong operator topology, by using the following facts
\begin{eqnarray*}
 R=\overline{\cup_{j} A_{2j}}^{\,\text{sot}}, \quad \quad  M^{(k+1)}_n \ot R=\overline{\cup_{j} M^{k+1}_n \ot A_{2j}}^{\,\text{sot}},
\end{eqnarray*}
we conclude that $\Gamma_k \circ \beta_k (R)=M^{(k+1)}_n \ot R$, which completes the proof.\qed
\end{prf}
\smallskip

\noindent\textbf{Proof of \Cref{diagonal subfactor}:} Consider the cyclic permutation $\sigma_1=\sum_{i=1}^{n-1}E_{i,i+1}+E_{n1}$. Let $\{\alpha_k\}$ be a sequence of inner automorphisms of the hyperfinite $\text{II}_1$ factor $R$ defined by  
\begin{eqnarray}\label{guru10}
\alpha_k(x) = \text{Ad}_{\sigma_1^{(k)}}(x)\,\,\mbox{ for all }x\in R,
\end{eqnarray}
where $\sigma^{(k)}_1$ is the $k$-times tensor product of the permutation matrix $\sigma_1$. By \Cref{outer automorphism}, the sequence $\{\alpha_k\}$ pointwise converges to some $\theta\in\text{Out}(R)$ with outer period $n$. For any $k\in\bbn\cup\{0\}$ and $\ell\in\bbn$, one observes that $\alpha_{k+\ell}(x)=\alpha_{k+1}(x)\,\forall\,x\in A_{2k}$. Therefore, it follows that $\theta(x)=\alpha_{k+1}(x)$ for all $x\in A_{2k}$. Now, consider the following projections
\[
e_1=\frac{1}{n} \sum_{i,j=1}^n E_{ij}\quad,\quad e_2 = \sum_{i=1}^nE_{ii} \otimes E_{ii}
\]
and define
\[
e_{2k+1}:=e_1 \otimes I_n^{(k)}\in A_{2k}\quad,\quad e_{2k}:=e_2 \otimes I_n^{(k-1)} \in A_{2k}.
\]
Since $\theta(x)=\alpha_{k+1}(x)$ for all $x\in A_{2k}$, in particular, using \Cref{guru10} we have the following\,:
 \begin{IEEEeqnarray}{lCl}\label{fix diagonal and projections}
     \theta(e_{2k+1}) &=& \alpha_{k+1}(e_{2k+1}) = e_{2k+1}\nonumber\\
     \theta(e_{2k}) &=& \alpha_{k+1}(e_{2k})=e_{2k}\nonumber\\ 
     \theta(D) &=& \alpha_1(D)=\text{Ad}_{\sigma_1}(D)\qquad\text{for all }\,D \in \dt.
 \end{IEEEeqnarray}
Now, recall the isomorphism $\Gamma_k\circ\beta_k$ from \Cref{make it diagonal}.
\smallskip

\noindent\textbf{Claim\,:} $\Gamma_k\circ\beta_k(R^{F_n, F_n}_k)=\mathscr{R}_k$ for each $k\in\bbn\cup\{0\}$.
\smallskip

\noindent Recall from \Cref{fourier algebra},
\[
N^{F_n, F_n}_{2k}=\text{Ad}_{(F_n)_{2k+1}(I_n \otimes (F_n)_{2k})}(I_n \otimes \dt \otimes I_n^{(k)})
\]
which is generated by the element $z:=\text{Ad}_{(F_n)_{2k+1}(I_n \ot (F_n)_{2k})}(I_n \ot \mathscr{D}_1 \ot I^{(k)}_{n})$. From \Cref{diagonal von Neumann algebra}, it follows  that
\[
z=\displaystyle{ \sum_{i, i_1,\cdots,i_k=1}^{n}\text{Ad}_{\sigma_{r_i}} (\mathscr{D}_1) \ot E_{i_1i_1}  \ot \cdots \ot E_{i_{k}i_{k}} \ot E_{ii}}.
\]
where $r_i\in \{0,1,2, \cdots,n-1\}$. Therefore, from \Cref{guru20}, we have the following\,:
\begin{IEEEeqnarray}{lCl}
\Gamma_k \circ \beta_k(z)&=&\displaystyle{\sum_{i, i_1,\cdots,i_k=1}^{n}  E_{i_1i_1}  \ot \cdots \ot E_{i_{k}i_{k}}\ot E_{ii} \ot \text{Ad}_{\sigma_{r_i}} (\mathscr{D}_1)}\qquad(\text{by \Cref{fix diagonal and projections}})\nonumber\\
&=&\displaystyle{ \sum_{i, i_1,\cdots,i_k=1}^{n}  E_{i_1i_1}  \ot \cdots \ot E_{i_{k}i_{k}}\ot E_{ii}\ot\theta^{r_i}(\mathscr{D}_1)}
\end{IEEEeqnarray}
(see \Cref{notation3}). Hence, we conclude that $\Gamma_k \circ \beta_k (z)\in\mathscr{R}_k$, and therefore,
\begin{equation}\label{send to diagonal algebra}
\Gamma_k \circ \beta_k (N^{F_n, F_n}_k) \subset \mathscr{R}_k.
\end{equation}
\noindent Similarly using \Cref{fix diagonal and projections}, one can verify that
\begin{IEEEeqnarray}{lCl}\label{projection1 maps to diagonal}
& & \Gamma_k \circ \beta_k(e_{2(k+1)+2j})=I^{(k+1)}_n \ot \theta^{r_i}(e_{2j}) \in \mathscr{R}_k,\nonumber\\
& & \Gamma_k \circ \beta_k (e_{2k+2j+1})=I^{(k+1)}_n \ot \theta^{r_i}(e_{2j-1}) \in \mathscr{R}_k.
\end{IEEEeqnarray}
From \Cref{send to diagonal algebra,projection1 maps to diagonal}, it follows that 
\[
\Gamma_k \circ \beta_k (\{N^{F_n, F_n}_k, \,e_{2(k+1)+i} \,: \, i\in \bbn\}) \subset \mathscr{R}_k.
\]
Now, using \Cref{vertex plus fourier} and the continuity of the isomorphism of $\Gamma_k \circ \beta_k$, we conclude that $\Gamma_k \circ \beta_k(R_{k}^{F_n, F_n}) \subset \mathscr{R}_k$. Therefore, we have the inclusion $\Gamma_k \circ \beta_k(R_k^{F_n,F_n}) \subset \mathscr{A}_k \subset M^{(k+1)}_n \ot R$ of type $II_1$ factors. By the multiplicativity of the Jones' index, we obtain $[\mathscr{A}_k :\, \Gamma_k \circ \beta_k(R_{k}^{F_n,F_n})]=1$,
which finishes the proof of the claim, and consequently, the first part of the theorem.

Finally, we have the following,
\begin{IEEEeqnarray}{lCl}\label{at end}
     (R^{F_n, F_n}_k)^{'} \cap R &=&(\Gamma_k \circ \beta_k)^{-1}\big(\mathscr{R}_k^\prime\cap(M^{(k+1)}_n \ot R)\big)\nonumber\\
     &=&(\Gamma_k \circ \beta_k)^{-1}(\Delta_{n^{k+1}} \ot I_n )\nonumber\\
     &=& \Delta_{n^{k+1}}
 \end{IEEEeqnarray}
which completes the proof.\qed

\noindent\textbf{Final discussion:} In \Cref{diagonal subfactor}, we have proved that when $u\in[F_n]$ and $v=F_n$, all the resulting vertex model subfactors are diagonal. It turns out that the moment we go beyond this situation, the outcome of the ``diagonal subfactor'' fails. Let us see counterexamples.

{\em The case of $u\notin[F_n]$}: Consider $u$ such that the spin model subfactor $R_u\subset R$ is of infinite depth. There are examples of such a situation, see \cite{Bur} for instance. In this case, $u\notin[F_n]$, as otherwise, $R_u\subset R$ would be of depth $2$ being a crossed product. Then, all the vertex model subfactors $R_k^{u,v}\subset R$, where $v$ is arbitrary, can not be diagonal as they all are of infinite depth, since $R_u\subset R$ is an intermediate subfactor with infinite depth.

{\em The case of $v\neq F_n$}: We have two subcases, namely, $u\in[F_n]$ and $u\notin[F_n]$. The situation of $u\notin[F_n]$ is already discussed above. So we are essentially left with the situation of $v\neq F_n$ and $u\in[F_n]$. In this regard, recall \Cref{to be used next}. There, $u=F_n,\,v=PF_n$ and $\big(R_0^{F_n,PF_n}\big)^\prime\cap R=\bbc$. This shows that $R_0^{F_n,PF_n}\subset R$ cannot be diagonal. Thus, we also have a counterexample for the case of $v\neq F_n$ and $u\in[F_n]$.
\smallskip

We conclude the article with the following natural question.
\smallskip

\noindent\textbf{Open question 4\,:} Consider a pair $(u,v)$ of complex Hadamard matrices of the same order and obtain the vertex model subfactors $R_k^{u,v}\subset R,\,k\in\bbn\cup\{0\}$, constructed in \Cref{theorem2}. Suppose that all these vertex model subfactors are diagonal. Is it necessary that $v=F_n$ and $u\in[F_n]$?
\smallskip

Observe that the `sufficiency' part is established in \Cref{diagonal subfactor}.


\newsection{Appendix}\label{Sec7}

Here, we provide a proof of Question B in \Cref{Sec4.1}, as promised there.

Let $u$ and $w$ be complex Hadamard matrices of order $n$. We show that for any $k, m\in \bbn$, the following subfactors
\[
R^{u,w}_{(m+2)k-1}  \subset  R^{u,w}_{(m+1)k-1} \subset R^{u,w}_{mk-1} \quad  \text{and} \quad R^{u,w}_{2k-1} \subset R^{u,w}_{k-1} \subset R 
\] 
in the nested sequence of vertex model subfactors  
\begin{equation*}
\cdots  \subset R^{u,w}_{(m+2)k-1}  \subset  R^{u,w}_{(m+1)k-1} \subset R^{u,w}_{mk-1} \subset \cdots \subset R^{u,w}_{2k-1} \subset R^{u,w}_{k-1} \subset R 
\end{equation*} 
need not be a downward basic construction. From \Cref{construction middle vertex model}, observe that the following quadruple 
\[
\begin{matrix}
 \bbc \ot M^{(m+1)k}_n  &\subset & R \cr \cup &\ &\cup\cr \bbc
&\subset{} &  R^{u,w}_{(m+1)k-1}
\end{matrix} 
\]
is a commuting square. Hence, it follows that $E^{R}_{R^{u,w}_{mk-1}}(\bbc \ot M^{(m+1)k}_n)=\bbc$. In particular, for the case of $u=w= F_n$, we have
\begin{eqnarray}\label{image of ce}
E^{R}_{R^{F_n,F_n}_{(m+1)k-1}}(\bbc \ot M^{(m+1)k}_n)=\bbc.
\end{eqnarray}
Now, from \Cref{construction middle vertex model} we have $N^{F_n, F_n}_{2(mk-1)} \subseteq \dt \ot M^{(mk)}_n \subseteq  R^{F_n,F_n}_{mk-1}$. Therefore, from \Cref{image of ce} we conclude that  
\begin{eqnarray}\label{image of cf}
 E^{R^{F_n,F_n}_{mk-1}}_{R^{F_n,F_n}_{(m+1)k-1}}(N^{F_n,F_n}_{2(mk-1))} )=\bbc.
\end{eqnarray}
To conclude the proof, it is now sufficient to construct an element $x \in R^{F_n, F_n}_{mk-1}$ such that
\[
f x f \neq E^{R^{F_n,F_n}_{mk-1}}_{R^{F_n,F_n}_{(m+1)k-1}}(x) f
\]
for any non-zero projection $f \in (R^{F_n, F_n}_{(m+1)k-1)})^{'} \cap R$. Consider 
\[
x=\text{Ad}_{(F_n)_{2mk-1}}\text{Ad}_{(I_n \otimes (F_n)_{2mk-2})}(I_n \otimes \mathscr{D}_1 \otimes I_n^{(mk)}) \in N^{F_n,F_n}_{2(mk-1)} \subseteq R^{F_n, F_n}_{mk-1}.
\]
From \Cref{diagonal von Neumann algebra}, it follows that $x \in  N^{F_n,F_n}_{2(mk-1)}$ is a diagonal unitary matrix where diagonal entries are from the set $\{1,\omega, \omega^2, \dots, w^{n-1}\}$ and $tr(x)=0$. Here, $\omega$ is a primitive $n$-th root of unity, and  $tr$ is the unique trace on the hyperfinite type $II_1$ factor $R$. Since $x \in  N^{F_n,F_n}_{2(mk-1)}$, using  \Cref{image of cf} we conclude that
\begin{eqnarray}\label{last}
E^{R^{F_n,F_n}_{mk-1}}_{R^{F_n,F_n}_{(m+1)k-1}}(x)=tr(x)=0.
\end{eqnarray}
Now, from \Cref{at end} we have
\begin{eqnarray*}
\big(R^{F_n, F_n}_{(m+1)k-1}\big)^\prime\cap R &= &\Delta_{n^{(m+1)k}}.
\end{eqnarray*}
Therefore, for any non-zero diagonal projection $f \in (R^{F_n, F_n}_{(m+1)k-1)})^{'} \cap R$, since $x \in  N^{F_n,F_n}_{2((m+1)k-1)}$ is a diagonal unitary matrix, it follows that
\begin{eqnarray*}
f x f&=&xf \not = 0.
\end{eqnarray*}
Observe the following from \Cref{last}
\[
E^{R^{F_n,F_n}_{mk-1}}_{R^{F_n,F_n}_{(m+1)k-1}}(x)f=tr(x)f=0.
\]
Thus, we conclude that
\[
f x f \neq   E^{R^{F_n,F_n}_{mk-1}}_{R^{F_n,F_n}_{(m+1)k-1}}(x)f\quad\forall\,\,0\neq f\in (R^{F_n, F_n}_{(m+1)k-1})^\prime\cap R,
\]
which concludes the proof.


\bigskip

\bigskip

\noindent {\em Department of Mathematics and Statistics},\\
{\em Indian Institute of Technology Kanpur},\\
{\em Uttar Pradesh $208016$, India}
\medskip

\noindent {\bf Keshab Chandra Bakshi:} keshab@iitk.ac.in, bakshi209@gmail.com\\
{\bf Satyajit Guin:} sguin@iitk.ac.in\\
{\bf Guruprasad:} guruprasad8909@gmail.com

\end{document}